\documentclass[10pt]{amsart}
\usepackage[all]{xy}
\usepackage{amsmath}
\usepackage{amsthm}
\usepackage{amsfonts}
\usepackage{amssymb}
\usepackage{amscd}
\usepackage{verbatim}
\newtheorem{theorem}{Theorem} [section]
\newtheorem{lemma} [theorem] {Lemma}
\newtheorem{proposition} [theorem] {Proposition}

\newtheorem{corollary} [theorem] {Corollary}

\begin{document}
\title{$L^2$-$\overline\partial$-cohomology groups
of some singular complex spaces }
\author{Nils  \O vrelid and Sophia Vassiliadou}
\thanks{The research of the second author is partially supported by NSF grant  DMS-0712795}
\address{Dept. of Mathematics\\University of Oslo\\
P.B 1053 Blindern, Oslo, N-0316 NORWAY}
\address{Dept. of Mathematics\\Georgetown University\\
Washington, DC 20057 USA} \email{nilsov@math.uio.no,\;
sv46@georgetown.edu}
\date{\today}
\begin{abstract} \noindent Let $X$ be a pure
$n$-dimensional (where $n\ge 2$) complex analytic subset in
$\mathbb{C}^N$ with an isolated singularity at $0$. In this paper
we express the $L^2$-$(0,q)$-$\overline\partial$-cohomology groups
for all $q$ with $1\le q\le n$ of a sufficiently small deleted
neighborhood of the singular point in terms of resolution data. We also obtain identifications 
of the $L^2$-$(0,q)$-$\overline\partial$-cohomology groups of the smooth points of $X$, in terms of resolution data,  when $X$ is either compact or an open relatively compact complex analytic subset of a reduced complex space with 
finitely many isolated singularities. 
\end{abstract}

\maketitle
\medskip
\noindent \section{Introduction}

\medskip
\noindent\noindent Let $X$ be a reduced pure $n$-dimensional complex
analytic set in $\mathbb{C}^N$ with an  isolated singularity at $0$
and let $X'$ denote the set of smooth points of
$X$. Let $(z_1,\cdots, z_N)$ be the coordinates in $\mathbb{C}^N,$ and
set $\|z\|:=(\sum_{j=1}^N |z_j|^2)^{\frac{1}{2}}$.  The set of smooth points $X'$ inherits a K\"ahler
metric from its embedding in $\mathbb{C}^N$, which we call the
ambient metric. Due to the incompleteness of the metric there are
many possible closed $L^2$-extensions of the
$\overline\partial$-operator originally acting on smooth forms on
$X'$. We consider the maximal (distributional)
$\overline\partial_{{\rm max}}$-operator.  For  positive $r$ we let $B_r:=\{z\in
\mathbb{C}^N;\;\|z\|<r \},\, X_r:=X\cap B_r,$ and  $X'_r:=X'\cap B_r$.  We shall choose an $R>0$
small enough, so that $bB_r$ intersects $X$ transversally  for all
$0<r< R$. Unless otherwise noted in what follows  by
$\overline\partial$ we shall mean $\overline\partial_{{\rm max}}$.
We define the local  (resp. global) $L^2$-$\overline\partial$- cohomology groups

$$
{H}^{p,q}_{(2)}(X'_r):= \frac{ \text{ker}\;(\overline\partial)\;\cap
L^{p,\,q}_{(2)}(X'_r)}{\text{Im}\;(\overline\partial)\cap
L^{p,\,q}_{(2)}(X'_r)},
$$

\medskip
\noindent (resp. ${H}^{p,q}_{(2)}(X'):= \frac{ \text{ker}\;(\overline\partial)\;\cap
L^{p,\,q}_{(2)}(X') }{\text{Im}\;(\overline\partial)\cap
L^{p,\,q}_{(2)}(X')}$
). 

\medskip
\noindent In \cite{FOV} we showed that the above local
$L^2$-$\overline\partial$-cohomology groups are finite dimensional
when $p+q<n$ and $q>0$ and zero when $p+q>n$. The idea of the proof
in the case $p+q<n$, was based on constructing complete K\"ahler
metrics to obtain a weighted $L^2$-solution for square-integrable,
$\overline\partial$-closed forms on $X'_r$, with compact support on
$X_r$ and identifying the obstructions to solving $\overline\partial
u=f$ on $X'_r$ to certain $L^2$-$\overline\partial$-cohomology
groups of ``spherical shells''  around $0$. Sharp regularity results
for $\overline\partial$ (which could yield  finite dimensionality
results for the above cohomology groups when $p+q\le n-2, \,q>0$)
have been obtained by Pardon and Stern for projective varieties with
isolated singularities in \cite{PS2}. We also presented in
\cite{FOV} various sufficient conditions on the
complex analytic set to guarantee that the local
$L^2$-$\overline\partial$-cohomology groups vanish. Our results were
most complete when $0$ was an isolated singular point in a
hypersurface $X$ and when $p+q\le n-1, \, 1\le q\le n-2$ ($n\ge 3$).
 In \cite{OV1} \, we proved finite dimensionality of
${H}_{(2)}^{n-1,1}(X'_r)$ using a global finite dimensionality
result of $L^2$-$\overline\partial$-cohomology groups on projective
varieties with arbitrary singularities.

\medskip
\noindent All of the results in \cite{FOV} were obtained while
working on the original singular space. This paper started as an
attempt to provide a short proof of the finite dimensionality of
$L^2$-Dolbeault cohomology groups of complex spaces with isolated
singularities by passing to an appropriate desingularization of $X$.
The second author had presented such results in conference talks
since 2006. Since then, new techniques have evolved to describe the
$L^{2}$-$({0,q})$-$\overline\partial$-cohomology groups in some
cases (see the work of Ruppenthal \cite{Rup1, Rup2} that
deals with cones over smooth projective varieties and his most recent preprint \cite{Rup3}). 
Using results from earlier papers of ours, some
classical theorems from algebraic geometry and singularity theory and some key observations from \cite{Rup3} and 
\cite{PS1}, we were able to obtain a rather complete
description of both local and global (the latter result when $X$ is compact or open relatively compact complex analytic set  in a reduced complex space with finitely many isolated singularities in $X$) $L^2$-$(0,q)$-$\overline\partial$-cohomology groups on $X'_r$ or $X'$ in terms of resolution data.  Earlier work  of Pardon (section 4 in \cite{P}) indicated the importance of such descriptions in 
understanding birational invariants of singular projective varieties. 

\medskip
\noindent The first  main result in the paper is the following
theorem:

\begin{theorem} Let $X$ be a complex analytic subset of\; $\mathbb{C}^N$ of pure
dimension $n\ge 2$ with an isolated singularity at $0$, and let
$\pi: \tilde{X}\to X$ be a desingularization. Then, there exists a
well-defined, linear mapping $\phi_*: H^q
(\tilde{X}_r,\,\mathcal{O})\to H^{0,q}_{(2)}(X'_r)$ such that
$\phi_*$ is bijective if\; $1\le q\le n-2$ and injective if\;
$q=n-1$.

\noindent Here $\tilde{X}_r=\pi^{-1}(X_r)$ and $X_r:=X\cap \{z\in
\mathbb{C}^N;\;\;\|z\|<r\}$.
\end{theorem}

\medskip
\noindent The above theorem generalizes results of Ruppenthal in
\cite{Rup2}. In that paper, he considered affine cones over smooth
projective varieties. For these varieties the exceptional locus of a
desingularization is a smooth submanifold of $\widetilde{X}$. We
impose no such restriction on the exceptional locus of the
desingularization. 
Key ingredient in the proof of Theorem 1.1 is a theorem of Stephen
Yau and Ulrich Karras (\cite{Yau2},\,\cite{Kar}) that describes the
local cohomology along exceptional sets. For complex analytic
subsets of $\mathbb{C}^N$ with an isolated singular point, the
exceptional locus of a desingularization is an exceptional set in
the sense of Grauert (see part $\alpha)$ in the  Characterization of
exceptional sets in section 3.1).

\medskip
\noindent The cokernel of the map $\phi_*$ will play a prominent
role in the paper. As we mentioned earlier, due to
the incompleteness of the metric, there are many $L^2$-extensions of
the $\overline\partial$-operator  acting on smooth forms on $X'_r$.
So far we have been considering the maximal (distributional)
extension. We can also consider the $L^2$-closure of
$\overline\partial$ acting on forms  with coefficients in
$C^{\infty}_0(\overline{X}_r\setminus \{0\})$. Let us denote this
extension by $\overline\partial^1$. We shall see in section 4, that
the cokernel of $\phi_*$ (or more precisely the dual of it) measures
somehow the obstructions to having $\overline\partial_{{\rm
max}}=\overline\partial^1$ at the level of holomorphic
$(n,0)$-forms.

\medskip
\noindent In January of 2010, we became aware of a recent preprint of Ruppenthal that appeared 
at the Erwin Schr\"odinger Institute  preprint series. Its purpose was to
describe explicitly the $L^2$-$\overline\partial$-cohomology of
compact complex spaces in terms of resolution data and thus answer a
conjecture by MacPherson on the birational invariance of the
$L^2$-Euler characteristic of projective varieties. After having
seen his preprint and using lemma 6.2 from \cite{Rup3}, we were able
to strengthen Theorem 1.1. More precisely we show the following:

\begin{theorem} Let $X$ be a  complex analytic
subset of\; $\mathbb{C}^N$ of  pure  dimension $n\ge 2$ with an isolated
singularity at $0$. Let $\pi:\tilde{X}\to X$ be a
desingularization such that the exceptional locus $E$ of $\pi$ is
a simple, normal crossings divisor. Let
$Z=\pi^{-1}(\text{Sing}\,X)$ be the unreduced exceptional divisor
of the resolution, let the support of $Z$ be denoted by $|Z|:=E$
and let $D:=Z-|Z|$. Then, there exists a natural surjective linear map

$$T: H^{n-1}(\tilde{X}_r,\,\mathcal{O}(D))\to
H^{0,n-1}_{(2)}(X'_r)
$$

\noindent whose kernel is naturally isomorphic to $H^{n-1}_E
(\tilde{X}_r,\,\mathcal{O}(D))$. Here $H^{n-1}_E(\tilde{X}_r,\,\mathcal{O}(D))$ means 
cohomology with support on $E$.
\end{theorem}

%$$
%\text{dim}_{\mathbb{C}}
%H^{0,n-1}_{(2)}(X'_r)=\text{dim}_{\mathbb{C}}
%H^{n-1}(\tilde{X}_r,\,\mathcal{O}(D))-\text{dim}_{\mathbb{C}}
%H^1(\tilde{X}_r,\,\mathcal{O}(-D)\otimes \mathcal{K}_{\tilde{X}_r})
%$$

\noindent 
\medskip
\noindent As a corollary of theorem 1.2, we recover Theorem 7.1
from \cite{Rup3} (for $q=n-1$). This theorem asserts that when the
line bundle associated to the divisor $-D=|Z|-Z$ is locally
semi-positive with respect to $X$, then $H^{0,q}_{(2)}(X'_r)\cong
H^q(\tilde{X}_r, \, \mathcal{O}(D))$ for all $0\le q\le n$.
Indeed, using Serre duality and Takegoshi's twisted vanishing
theorem \footnotemark\footnotetext{For a proper, generically finite
to one holomorphic map $p: X\to Y$ where $X$ is a complex
connected manifold and $Y$ is a  reduced analytic space, this
relative vanishing theorem was already known to A. Silva (see A.2
Lemma in \cite{Sil}).} (Torsion freeness of the main theorem in
the introduction of \cite{Tak1}), we see that
$H^{n-1}_c(\tilde{X}_r,\,\mathcal{O}(D))=0$ in this case. A result
by Karras will guarantee the isomorphism between
$H^{n-1}_E(\tilde{X}_r,\,\mathcal{O}(D))$ and
$H^{n-1}_c(\tilde{X}_r,\,\mathcal{O}(D))$, which combined with
Theorem 1.2  will yield the desired isomorphism
$H^{0,n-1}_{(2)}(X'_r)\cong H^{n-1}(\tilde{X}_r, \,
\mathcal{O}(D))$. We can also recover Ruppenthal's result for all  
$q\le n-2$ (see Remark 4.5.1  in section 4). 

\medskip
\noindent In order to prove Theorem 1.2, we construct a
non-degenerate pairing

\begin{eqnarray} \label{eq:1dp}
\dfrac{H^{0,n-1}_{(2)}(X'_r)}{\phi_*(H^{0,n-1}(\tilde{X}_r))} \times
\dfrac{\text{kern}(\overline\partial)^{n,0}}{\text{kern}(\overline\partial^1)^{n,0}}\to
\mathbb{C}
\end{eqnarray}

\medskip
\noindent where $\overline\partial^1$ is as above. In \cite{FOV} we showed that the map $j_*:
H^{0,\,n-1}_{(2)}(X'_r)\to H^{n-1}(X'_r,\,\mathcal{O})$ induced by
the natural inclusion $j: L^{0,\,n-1}_{(2)}(X'_r)\to
L^{0,n-1}_{2,\,\text{loc}}(X'_r)$ is injective. An understanding of
the $\text{Im}\,j_*$ will turn out to be instrumental in the
construction of the map $T$. We will therefore present some
necessary and sufficient conditions to describe elements in
$\text{Im}\,j_*$ (using Lemma 6.2 in \cite{Rup3} and $(\ref{eq:1dp})$).
Now, there exists a natural map $\ell_*:H^{n-1}(\tilde{X}_r,\,\mathcal{O}(D))\to
H^{n-1}(X'_r,\,\mathcal{O})$. Using a twisted version of an $L^2$-Cauchy problem we will show that
$\text{Im}\,j_*\subset \text{Im}\,\ell_*$ and construct a map $S: H^{0,n-1}_{(2)}(X'_r)\to 
H^{n-1}(\tilde{X}_r, \mathcal{O}(D))$. Then the proof of Theorem 
1.2  will be based on the following key observation: the map $\ell_*$ is surjective
on the $\text{Im}\,j_*$. The composition ${j_*}^{-1}\circ \ell_*$
will be the desired map $T$ and $T\circ S=Id$. 

\medskip
\noindent When $q=n$, we can easily show that the map $\phi^n_*: H^{n}(\tilde{X}_r, \,\mathcal{O}) \to H^{0,n}_{(2)}(X'_r)$ described by $\phi^n_*([g])=[(\pi^{-1})^* g]$ is surjective. Since $\tilde{X}_r$ contains no compact $n$-dimensional irreducible components,  by  Siu's theorem (\cite{Siu})  we have  $H^{n}(\tilde{X}_r, \,\mathcal{O})=0$. Hence, $H^{0,n}_{(2)}(X'_r)=0$. 

\medskip
\noindent With a little bit more work, we can obtain global versions of Theorems 1.1 and 1.2. More precisely, let 
$X$ be a  pure $n$-dimensional, relatively compact domain  in a reduced complex analytic space $Y$. We give $\text{Reg}\,Y$ a hermitian metric compatible with local embeddings. Assume that $\overline{X}\cap \text{Sing}\,Y=:\Sigma=\{a_1,a_2,\cdots,a_m\}\subset X$. Let $\pi:\tilde{Y}\to Y$ be a desingularization such that $E=\pi^{-1}(\Sigma)$ is a normal crossings reduced  divisor  in $\tilde{X}=\pi^{-1}(X)$. Let $Z:=\pi^{-1}(\Sigma)$ be the unreduced exceptional divisor and $D:=Z-E$. Give $\tilde{Y}$ a non-degenerate hermitian metric. Let $H^{0,q}_{(2)}(\tilde{X},\,\mathcal{O}(D))$ denote the $L^2$-cohomology of $(0,q)$-forms in $\tilde{X}$ with values in $L_D$, the holomorphic line bundle associated to the divisor $D$ (see Remark 2.2.2 in section 2).  Then we have

\begin{theorem} The map $\phi_*: H^{0,q}_{(2)}(\tilde{X})\to H^{0,q}_{(2)}(X')$, defined by $\phi_*([f])=[(\pi^{-1})^*\,f]$ is an 
isomorphism, when $1\le q\le n-2$  and  where $X':=X\setminus \Sigma$.
\end{theorem} 

\medskip
\noindent 
\begin{theorem}There exists a natural surjective map $\tilde{T}: H^{0,\,n-1}_{(2)}(\tilde{X},\,\mathcal{O}(D))\to H^{0,n-1}_{(2)}(X')$,
whose kernel is naturally isomorphic to $H^{n-1}_E(\tilde{X},\,\mathcal{O}(D))$ and where $X':=X\setminus \Sigma$.  
\end{theorem}

\medskip
\noindent Let us point out that in  the most interesting cases, i.e. when  $X$ is compact or $\partial X$ is smooth,   strongly pseudoconvex submanifold of $\text{Reg}\,Y$, we have  $H^{0,q}_{(2)}(\tilde{X},\,\mathcal{O}(F))\cong H^{q}(\tilde{X},\,\mathcal{O}(F))$ for $q>0$ and $F$ any holomorphic line bundle,  so Theorems 1.1 and 1.2 carry over verbatim. In \cite{Rup3}, Ruppenthal proved (Theorem 1.6) that when the line bundle associated to the divisor $-D$ is locally semi-positive with respect to $X$, then $H^{0,n-1}_{(2)}(X')\cong H^{n-1}(\tilde{X},\,\mathcal{O}(D))$
\footnotemark\footnotetext{More precisely, Theorem 1.6 in \cite{Rup3}, states that when the line bundle associated to the divisor $-D$ is locally semi-positive with respect to $X$, then for all $q$, with $0\le q\le n$ one has 
$H^{0,q}_{(2)}(X')\cong H^q(\tilde{X},\,\mathcal{O}(D))$.}.  This follows from Theorem 1.4 taking into account Takegoshi's or Silva's 
relative vanishing theorem and Karras' results. For projective surfaces with isolated singularities, we can say more: 

\begin{corollary}  Let $X$ be a projective  surface  with finitely many isolated singularities. Then the map $\tilde{T}: H^{0,1}_{(2)}(\tilde{X},\mathcal{O}(D))\to H^{0,1}_{(2)}(X')$ is an isomorphism (the right-hand side $L^2$-cohomology is computed with respect to the restriction of the Fubini-Study metric in $X'$). 
\end{corollary}

\medskip
\noindent This Corollary was first conjectured by Pardon in \cite{P}, while studying MacPherson's conjecture. It appeared later as a special 
case of  Theorem B in \cite{PS1}.  A key observation from the Appendix in \cite{PS1} along with Theorem 1.4 will help us settle Pardon's conjecture  in the case of projective surfaces with isolated singularities and bypass the difficulties that were encountered with the proof of Theorem B in \cite{PS1}.  It would be interesting to determine whether the kernel of $\tilde{T}$ vanishes for higher dimensional projective varieties with an isolated singularity (Professor Koll\'ar offered some insight on when this vanishing could occur; see  Remark 5.2.3   in section 5).  In that case the global cohomology group $H^{0,n-1}_{(2)}(X')$ would be isomorphic to $H^{n-1}(\tilde{X},\,\mathcal{O}(D))$. Correspondingly, this $L^2$-Dolbeault cohomology group would  not be a birational invariant. 

\medskip
\noindent  Now we follow the assumptions and notation as in the paragraph just above Theorems 1.3, 1.4 and consider the case where $X$ is compact or $\partial X$ is smooth strongly pseudoconvex submanifold of $\text{Reg}\,Y$.  The  map $\phi^n_*: H^{n}(\tilde{X},\,\mathcal{O})\to H^{0,n}_{(2)}(X')$  defined by $\phi^n_*([f])=[\,(\pi^{-1})^*\,f\,]$  is easily seen to be surjective. Let $i^n_*: H^{n}(\tilde{X},\,\mathcal{O})\to H^n(\tilde{X}, \,\mathcal{O}(D))\cong H^{0,n}_{(2)}(\tilde{X},\,L_D)$ be the map on cohomology induced by the sheaf inclusion  $i: \mathcal{O}\to \mathcal{O}(D)$.  We will  show 

\begin{corollary} With $X,\tilde{X}, D, \phi^n_*,\,i^n_*$ as above we have $\text{kern}(\phi^n_*)=\text{kern}(i^n_*)$ and $H^{0,n}_{(2)}(X')\cong H^n(\tilde{X},\,\mathcal{O}(D))$.
\end{corollary} 
 
\medskip
\noindent  The kernel of $i^n_*$ can be computed using standard long exact sequences on cohomology and cohomology with support on $E$.  Thus one of the benefits of the above Corollary is that it allows us to describe  the kernel of $\phi^n_*$ which in some sense measures  the difference between the $L^{0,n}_{(2)}$-$\overline\partial_{min}$-cohomology group on $X'$  (which is isomorphic to $H^n(\tilde{X},\,\mathcal{O})$),  and the corresponding cohomology group using the $\overline\partial_{max}$-operator (i.e. $H^{0,n}_{(2)}(X')$).  
\medskip  

\medskip
\noindent The organization of the paper is as follows: Apart for
some preliminaries, in section 2 we will  give a short proof of the
finite dimensionality of $L^{p,q}_{(2)}$-$\overline\partial$-cohomology groups of small
deleted neighborhoods of the singular point $0$. In section 3 we
prove Theorem 1.1. Section 4 contains the proof of Theorem 1.2. Section 5 contains the proofs of the global theorems and section 6, the identification of $H^{0,n}_{(2)}(X')$ with $H^n(\tilde{X}, \mathcal{O}(D))$.  In section 7 we discuss the vanishing or not of some local $L^2$-$\overline\partial$-cohomology groups of some complex spaces $X$ with isolated singularities.

\medskip
\noindent{\bf{Acknowledgements:}} The first author wants to thank Jean Ruppenthal for a clarifying e-mail exchange about the results in \cite{Rup3}. This paper started  while
the second author was on sabbatical leave from Georgetown University in the Spring of 2009.  Part of the work was done while she was visiting the Department of Mathematics at the  University of Chicago, the Institute for Advanced Studies in Princeton and the Institute of Mathematics at the University of Oslo. She would like to thank these institutions and in particular the several complex variables group at the University of Oslo for their hospitality. She would also like to thank Mark de Cataldo and Mihnea Popa for helpful correspondence, Professor J\'anos Koll\'ar for numerous suggestions in determining  $H^{n-1}_E(\tilde{X},\,\mathcal{O}(D)),$  Tom Graber for helpful remarks on the cohomology of projectivized vector bundles  and Tom Haines for fruitful discussions. She acknowledges financial support from IAS under NSF grant DMS--0635607. Any opinions, findings and conclusions or recommendations expressed in this material are those of the author and do not necessarily reflect the views of the National Science Foundation. 

\section{Preliminaries}

\subsection{Desingularization and pull-back metrics} Our results in \cite{FOV} were obtained while working mostly on the
original singular space. However, we can desingularize $X$, i.e.
consider a proper, holomorphic, surjective map $\pi:\tilde{X}\to X$
such that $\tilde{X}$ is smooth, $\pi:\tilde{X}\setminus E\to
X\setminus {\rm Sing}X$ is a biholomorphism and $E=\pi^{-1}({\rm
Sing X})$ is a divisor with normal crossings (we only  need this extra condition on the exceptional locus for the proof
of Theorems 1.2 and 1.4).  Since the singular
locus of $X$ consists of one point we can cover $E$ by finitely many
coordinate charts $(U_i,z)$ with $i=1,\cdots, M$ and near each
$x_0\in E$ we can find local holomorphic coordinates
$(z_1,\dots,z_n)$ in terms of which $E\cap U_i$ is given by
$h_i(z)=z_1\cdots z_{n_i}=0,$ where $1 \leq n_i \leq n$.

\medskip
\noindent Let $\sigma$ be a positive definite metric on $\tilde{X}$.
We can then consider, volume element $d \tilde{V}_{\sigma}$ and
pointwise norms/norms on $\Lambda^{\cdot} T \tilde{X}$ and
$\Lambda^{\cdot} T^*\tilde{X}.$
\medskip
\noindent For every open subset $U$ of $\tilde{X}$, let
\;$\mathcal{L}^{p,q}(U)$ be

$$
\mathcal{L}^{p,q}(U):=\{u\in L^{p,\,q}_{2,\,{\rm loc}}(U);\;\;
 \overline{\partial}u\in  L^{p,\,q+1}_{2,\,{\rm loc}}(U)\}
$$

\noindent  and for each open subset $V\subset U$, let $r^{U}_{V}:
\mathcal{L}^{p,q}(U)\to \mathcal{L}^{0,q}(V)$ be the obvious
restriction maps. Here square-integrability is with respect to the
metric $\sigma$. Then the map $u\to \overline\partial u$ defines an
$\mathcal{O}_{\tilde{X}}$-homomorphism $\overline\partial:
\mathcal{L}^{p,q}\to \mathcal{L}^{p,q+1}$ and the sequence

$$
0 \to \mathcal{O}_{\tilde{X}}\to \mathcal{L}^{p,0}\to
\mathcal{L}^{p,1}\to \cdots \to \mathcal{L}^{p,n}\to 0
$$

\noindent is exact by the local Poincar\'e lemma for
$\overline\partial$. Since each $\mathcal{L}^{p,q}$ is closed under
multiplication by smooth cut-off functions we have a fine resolution
of $\mathcal{O}_{\tilde{X}}$.

\smallskip
\noindent We introduce some notational convention: For the manifold
$\tilde{X}$, $\gamma$ will always denote a positive semi-definite
hermitian metric on $\tilde{X}$, which is generically definite. More
specifically in this paper we shall let $\gamma$ denote the
pull-back of the ambient metric on $X'$. It degenerates along a
divisor $D_{\gamma}$ supported on the exceptional divisor $E$. One
is faced with the hard task of understanding how the pull back of
the ambient metric looks like on $\tilde{X}$. This has been done by
Hsiang-Pati \cite{HP}, Nagase \cite{N} for projective surfaces with
isolated singularities and recently by Taalman \cite{T} (following
an idea of Pardon and Stern \cite{PS3}) for three-dimensional
projective varieties with isolated singularities. Youssin in
\cite{Y} considered desingularizations $(\tilde{X}, \pi)$ of $X$
that factor through the Nash blow-up of $X$ and found a way to
describe the pull-back of forms defined on $X'$ with measurable
coefficients and square-integrable with respect to the ambient
metric, in terms of data on $\tilde{X}$. Similar descriptions of
such forms for projective surfaces with isolated singularities
appeared in the 1997 preprint of Pardon and Stern \cite{PS3}.

\smallskip
\noindent

\medskip
\noindent \subsection{Locally free sheaves and twisted cohomology
groups} Let us consider any effective divisor $D=\sum_{i=1}^m\,d_i
\,E_i$ where $d_i\in \mathbb{N}$ and where $\{E_i\}_{i=1}^m$ are the
irreducible components of $E=\pi^{-1}(\text{Sing}\,X)$. By
$\mathcal{O}(D)$ we denote the sheaf of germs of meromorphic
functions $f$ such that $\text{div}\,(f)+D\ge 0$. If $\{U_a\}$ is a
covering of $\tilde{X}$ and $u_a$ is a meromorphic function on $U_a$
such that $\text{div}\,(u_a)=D$ on $U_a$, then
$\mathcal{O}(D)_{\upharpoonright \,U_a}=u_a^{-1}\,\mathcal{O}$.
Hence $\mathcal{O}(D)$ is a locally free sheaf of rank 1. This sheaf
can be identified with the sheaf of sections of a line bundle $L_D$ over $\tilde{X}$ defined by
the cocycle $g_{ab}:=\frac{u_a}{u_b}\in \mathcal{O}^*(U_a\cap U_b)$.
In fact there is a sheaf homomorphism $\mathcal{O}(D)\to
\mathcal{O}(L_D)$ defined by

$$
\mathcal{O}(D)(W)\ni f \mapsto s_f\in
\mathcal{O}(L_D)(W)\;\;\;\;\text{with}\;\;\theta_a(s_f)=f\,u_a
\;\;\;\text{on}\;\;W\cap U_a
$$

\noindent where $\theta_a$ is the corresponding trivialization of
${L_D}_{\upharpoonright U_a}$. The constant function $f=1$ induces a
meromorphic section $s$ of $L_D$ such that
$\text{div}\,(s)=\text{div}\,(u_a)=D$. Since $D\ge 0$, the section
$s$ is holomorphic and its zero set $s^{-1}(0)$ is  the support of
$D$, usually denoted by $|D|$.  Hence, we can identify  sections in $\Gamma(U,\,\mathcal{O}(D))$
with sections in $\Gamma(U,\,\mathcal{O}(L_D))$ via the isomorphism $f\to f\otimes s$. The inverse of this map is given by taking any section $ \Gamma(U,\,\mathcal{O}(L_D))\ni A \to A\;^{.}\,s^{-1}\in \Gamma(U,\,\mathcal{O}(D))$. Locally this map is described  by sending $A=u\otimes e\to s^{-1}(e)\,u,$  where $e$ is a local holomorphic frame for $\mathcal{O}(L_D)$ and $s^{-1}$
is a meromorphic section of $L_{-D}\cong L_D^*$ (the dual of $L_D$)  satisfying $s\cdot s^{-1}:=s^{-1}(s)=1$.

\medskip
\noindent For any open set $U\subset \tilde{X}$ we set

$$L^{p,q}_{2,\,\text{loc}}(U,\,\mathcal{O}(D)):=\{f\in
L^{p,q}_{2,\,\text{loc}}(U\setminus E)\;\;|\;\;\; \chi\,f\in
L^{p,q}_{2,\,\text{loc}}(V)\;\;\text{for\;\;all\;\;}\;\;V^{\text{open}}\subset
U\;\;\text{and}\;{\;\forall}\;\;\chi\in\mathcal{O}(-D)(V)\;\}.
$$ 

\medskip
\noindent{\bf{Remark 2.2.1}} In principle one could define $L^{p,q}_{2,\,\text{loc}}(U,\,\mathcal{O}(D))$ to consist of all forms $f\in L^{p,q}_{2,\,\text{loc}}(U\setminus |D|)$ such that $ \chi\,f\in
L^{p,q}_{2,\,\text{loc}}(V)\;\;\text{for\;\;all\;\;}\;\;V^{\text{open}}\subset
U\;\;\text{and}\;{\;\forall}\;\;\chi\in \mathcal{O}(-D)(V)$. But then for points 
$x\in E\setminus |D|$ one sees that $f$ extends as an $L^{p,q}_{2,\,\text{loc}}$ form across these points.  Hence, we do not lose any information by defining $L^{p,q}_{2,\,\text{loc}}(U,\,\mathcal{O}(D))$ the way we did before the Remark.  

\medskip
\noindent Similarly, for a sufficiently small, relatively compact open neighborhood $U$ of $E$ in $\tilde{X}$, one can define the following spaces

$$L^{p,q}_{(2)}(U,\,\mathcal{O}(D)):=\{f
\in L^{p,q}_{2,\,\text{loc}}(U\setminus E) \;\;|\;\;\; u_a\,f\in
L^{p,q}_{(2)}(U\cap U'_a)\;\;\;\text{for\;all\;a}\}
$$

\medskip\noindent where $\{U'_a\}$ is a finite open covering of $E$ and if $D=\sum d_j E_j$ and
$g_{j,a}$ is the local generator of the ideal sheaf of $E_j$ over $U'_a$, then $u_a:=\prod {g^{d_j}_{j,a}}.$\footnotemark\footnotetext{In section 5, we will consider $L^{p,q}_{(2)}(U,\,\mathcal{O}(D))$ for 
open sets $U\subset\subset \tilde{X}$  (or relatively compact in $\tilde{Y}$ where $Y,\,\tilde{Y}$ are as in the first paragraph above Theorem 1.3 in the introduction). Then the above definition can be reformulated by saying that $L^{p,q}_{(2)}(U,\,\mathcal{O}(D))$ consists of those $f\in L^{p,q}_{2,\,\text{loc}} (U\setminus E)$ such that $u_a\,f\in L^{p,q}_{(2)}(U\cap U'_a)$ for all $a$, and where $u_a$ is a generator of the ideal sheaf $\mathcal{O}(-D)$ of $D,$  on a neighborhood of $\overline{U'_a}$  for all $a,$ and $\{U'_a\}_{a\le m}$ is a covering of $\overline{U}$.}
In the above definition, square-integrability is with respect to any non-degenerate metric $\sigma$ on $\tilde{X}$. 
 It is clear from the definitions that for such a $U$ we have: $L^{p,q}_{(2)}(U,\,\mathcal{O}(D))\hookrightarrow L^{p,q}_{(2),\,\text{loc}}(U,\,\mathcal{O}(D))$. Using a partition of unity $\{\rho_a\}$ subordinate to the covering $\{U'_a\}$,  we can define a norm on this space: 

$$
\|f\|_{L^{p,q}_{(2)}(U,\,\mathcal{O}(D))}:=\left(\int_U \sum_a \rho_a |u_a\,f|^2\, dV\right)^{\frac{1}{2}}
$$

\noindent
This definition seems to depend on the covering $\{U'_a\}$, \, the partition of unity $\{\rho_a\}$, and the choice of the local defining function  for the divisor $D$. Since $\overline{U}$ is bounded, by passing to a slightly smaller covering of $U$, we will see that the corresponding norms, if we choose different coverings, defining functions for $D$ and partitions of unity, would be equivalent.

\medskip
\noindent Now the map $U\mapsto
\mathcal{L}^{p,q}(\mathcal{O}(D))(U):=\{f\in
L^{p,q}_{2,\,\text{loc}}(U,\,\mathcal{O}(D))\;\;\text{such\;\;that}\;\;\overline\partial\,f\in
L^{p,\,q+1}_{2,\,\text{loc}}(U,\,\mathcal{O}(D))\}$ (here
$\overline\partial$ is with respect to open subsets of
$\tilde{X}\setminus E$) is a fine sheaf on $\tilde{X}$ and

$$
0\to \Omega^p_{\tilde{X}}\otimes \mathcal{O}(D)\to
\mathcal{L}^{p,0}(\mathcal{O}(D))\overset{\overline\partial}\to
\mathcal{L}^{p,1}(\mathcal{O}(D))\overset{\overline\partial}\to
\cdots \overset{\overline\partial}\to
\mathcal{L}^{p,n}(\mathcal{O}(D))\to 0
$$

\noindent is a fine resolution of $\Omega^p_{\tilde{X}}\otimes
\mathcal{O}(D)$. To see this we can argue as follows: For $x\in
U_a$, the maps of germs $f_x\to (u_a\,f_x)\otimes {u_a}^{-1}$ from
${\mathcal{L}^{p,q}(\mathcal{O}(D))}_x\to
\mathcal{L}^{p,q}_x\otimes_{\mathcal{O}_x}\,\mathcal{O}(D)_x$ are
independent of $a$, where $\mathcal{L}^{p,q}$ are defined in section
2.1. These maps of germs define sheaf isomorphisms
$\mathcal{L}^{p,q}(\mathcal{O}(D))\to
\mathcal{L}^{p,q}\otimes_{\mathcal{O}} \mathcal{O}(D)$, commuting
with $\overline\partial$ and $\overline\partial\otimes Id$
respectively. Moreover, the operation
$-\otimes_{\mathcal{O}}\mathcal{O}(D)$ preserves exact sequences,
since $\mathcal{O}(D)$ is a locally free sheaf over $\mathcal{O}$.
Hence the cohomology of
$\left(\Gamma(\tilde{X}_r,\,\mathcal{L}^{p,\bullet}(\mathcal{O}(D))),\,\overline\partial\right)$
is $H^{\bullet}(\tilde{X}_r,\,\Omega^p\otimes \mathcal{O}(D))$ for
any $p\ge 0$. 

%\medskip\noindent We can also consider for each open subset $U$, the spaces 

%$$\mathcal{L}^{p,q}_{(2)}(U,\,\mathcal{O}(D)):=\{f\in L^{p,q}_{(2)}(U,\,\mathcal{O}(D))\;\;|\;\;\overline\partial f %\in L^{p,q+1}_{(2)}(U,\,\mathcal{O}(D))\}.$$ 

%\noindent Again, here $\overline\partial$ is with respect to open subsets of
%$\tilde{X}\setminus E$. Therefore one obtains a complex

%$$
%....\to \mathcal{L}^{p,q-1}_{(2)}(U,\,\mathcal{O}(D))\overset{\overline\partial}\to %\mathcal{L}^{p,q}_{(2)}(U,\,\mathcal{O}(D))\overset{\overline\partial}\to %\mathcal{L}^{p,q+1}_{(2)}(U,\,\mathcal{O}(D))\overset{\overline\partial}\to...
%$$

%\noindent and from it one can consider also the corresponding cohomology groups $H^{p,q}_{(2)}(U,\,\mathcal{O}(D))$.

\subsubsection{\bf{An alternative characterization of $L^{p,q}_{(2)}(U,\,\mathcal{O}(D))$}} 
In section 4 of this paper we would need another realization of $L^{p,q}_{(2)}(U,\;\mathcal{O}(D))$ for $U$ a smoothly bounded strongly pseudoconvex neighborhhod of $E$ in $\tilde{X}$. We would like to identify this space with 
the square-integrable  sections of $\wedge^{p,q}\,T^*_{\tilde{X}}\otimes L_D$ over $U$, where $L_D$ is the holomorphic line bundle associated to the divisor $D$.  We would also need in section 4, some general results about differential operators acting on sections of holomorphic line bundles, 
cohomology groups with coefficients in line bundles etc.  
In this section we will systematically discuss these notions. Let $\tilde{X}$
be given a non-degenerate metric $\sigma$ and let $F$ be a holomorphic line bundle endowed with a Hermitian metric $h$. Let $C^{\infty}_{p,q}(U,\,F):=C^{\infty}(U,\,\wedge^{p,q} T^* {\tilde{X}}\otimes F)$ denote the space of smooth $(p,q)$-forms in $U$  with coefficients in $F$, \; $C^{\infty}_{p,q}(\overline{U},\,F)$ denote the smooth up-to the boundary of $U$, $(p,q)$-forms with coefficients in $F$ and let $D^{p,q}(U,\,F)$ the compactly supported sections with coefficients in $F$. Using a trivialization $\theta_U: {F}_{\upharpoonright_{U}}\to U\times \mathbb{C}$ we can choose a frame $e(x):=\theta_U^{-1}(x,1)$ of $F$. Locally for each $x\in U$, any element $A\in C^{\infty}_{p,q}(U,F)$ \; can be written as $A=\phi\otimes e$ in a smaller neighborhood $W\subset U$ of $x$  where 
$\phi\in C^{\infty}_{p,q}(W)$ and $e\in \mathcal{O}(F)(W)$.  Let $\tau: F\to F^*$ be the conjugate-linear isomorphism of $F$ onto its dual $F^*$ defined by $\tau(e)(e'):=h(e',e)$ whenever $e,e'\in F_x$. The dual bundle $F^*$ is given the metric $h^*:=h^{-1}$ that makes $\tau$ an isometry. Then we can define the generalized Hodge-star-operator 

\begin{eqnarray}\label{eq:dfntinp}
\overline{*}_F\;&:& C^{\infty}_{p,q}(U,\,F)\to C^{\infty}_{n-p,\,n-q}(U,\,F^*) \notag\\ 
\newline \notag\\
&\;& \overline{*}_F(\phi\otimes e)=\overline{*}\,\phi \otimes \tau(e) \\ \notag
\end{eqnarray}

\noindent where $\phi\in \wedge^{p,q} T^*_x U$ and $e\in F_x$. 

\medskip
\noindent For sections
$A\in C^{\infty}_{p,q}(U,\,F)$ we can easily check that the following equality holds: $\overline{*}_{F^*}\,\overline{*}_F A=(-1)^{p+q}\,A$, where $\overline{*}_{F^*}$ is the Hodge-star operator associated to
$F^*$.

\medskip\noindent We can also define a wedge product $\wedge: C^{\infty}_{p,q}(U,\,F)\times C^{\infty}_{r,s}(U,F^*) \to C^{\infty}_{p+r,\,q+s}(U,\,\mathbb{C})$ described by 

\begin{equation}\label{eq:dfnwedge}
(\phi\otimes e)\wedge (\psi\otimes f)=\phi\wedge \psi\;f(e)
\end{equation}

\noindent where $A:=\phi\otimes e$ and $B:=\psi\otimes f$ are the local descriptions of two sections $A\in C^{\infty}_{p,q}(U,\,F)$ and $B\in C^{\infty}_{r,s}(U,F^*)$ and where $e,\,f$ are local frames for $F,\,F^*$ respectively.  

\medskip\noindent Using the metric $\sigma$ on $\tilde{X}$, the hermitian metric $h$ on $F$ and the local description of elements in $C^{\infty}_{p,q}(U,\,F)$ we can define a pointwise inner product for two elements $A,\,B\in C^{\infty}_{p,q}(U,\,F)$

\begin{equation}\label{eq:df1pointin}
<A,\,B>_{F,x}=h(e,\,e)\; <\phi,\,\psi>_{\sigma,\,x}\,
\end{equation}

\noindent where $A=\phi\otimes e$ and $B=\psi\otimes e$ in a small neighborhood $W\subset U$ of $x$ and $<\;,>_{\sigma,\,x}$ is the standard pointwise inner product on $\tilde{X}$ arising from the metric $\sigma$. By integrating with respect to the volume element $dV_{\sigma}$ we obtain a global $L^2$ inner product on $U$. 

\medskip\noindent  
For any two sections
$A,\,B\in C^{\infty}_{p,q}(U,\,F)$ given locally by $A=\phi\otimes e$ and $B=\psi\otimes e$ with $\phi,\,\psi$ smooth $(p,q)$-forms in smaller neighborhood of $x$ we have 

$$A\wedge \overline{*}_F\;B=\phi\wedge\overline{*}\psi\,\;(\tau(e))(e)=h(e,e)\,<\phi,\,\psi>_{\sigma,\,x}\,dV=<A,B>_{F,x}\,dV$$

\medskip\noindent As before, we obtain a global inner product on sections in $C^{\infty}_{p,q}(U,L_D)$ given by 

\begin{equation}\label{eq:dfninnp}
(A,\,B)_F=\int_U A\wedge \overline{*}_F B.
\end{equation} 

\medskip\noindent Let $\overline\partial_F=\overline\partial \otimes Id: C^{\infty}_{p,q}(U, F)\to C^{\infty}_{p,q+1}(U,\,F)$. Then we can define the formal adjoint 

$$\vartheta_F: C^{\infty}_{p,q}(U,F)\to C^{\infty}_{p,q-1}(U,\,F)$$

\noindent  via the identity $\vartheta_F:=-\overline{*}_{F^*}\,\overline\partial_{F^*}\,\overline{*}_F$, where by $\overline\partial_{F^*}$ we denote the $\overline\partial$  operator associated to the $F^*$. 

\medskip
\noindent Let $L^{p,q}_{(2)}(U,\,F)$ denote the completion of $D^{p,q}(U,\,F)$ under the inner product defined above. This completion is independent of the choice of the bundle metric $h$, with different choices of metrics leading to equivalent inner products.   The wedge product, inner product, the generalized Hodge $\overline{*}_F$ operator defined earlier for smooth sections, extend naturally to square-integrable sections. One also obtains various extensions of the operators $\overline\partial_F,\,\vartheta_F$ on $L^{p,q}_{(2)}(U,\,F)$  just as in the case of complex-valued forms. By abuse of notation we shall denote  the weak extension of $\overline\partial_F$ on $L^{\bullet,\,\bullet}_{(2)}(U,\,F)$ by $\overline\partial_F$ (instead of the cumbersome $(\overline\partial_F)_w$),    the minimal extension of $\overline\partial_F$ by  $\overline\partial_{F,\,min}$, the weak extension of $\vartheta_F$ on $L^{\bullet,\,\bullet}_{(2)}(U,\,F)$  by $\vartheta_{F,\,h}$ (instead of $(\vartheta_F)_w$)  and finally $\overline\partial^*_{F,\,h}$ will denote the Hilbert space adjoint of $\overline\partial_F$. Let 

$$
H^{p,q}_{(2)}(U,\,F):=\dfrac{\text{kern}(\overline\partial_F)\cap L^{p,q}_{(2)}(U,F)}{\text{Im}(\overline\partial_F)\cap L^{p,q}_{(2)}(U,\,F)}
$$

\noindent denote the $L^2$-cohomology groups with coefficients in $F$.

\medskip
\noindent {\bf{Remark 2.2.2}}  In sections 4 and 6 of the paper we will be considering  forms with coefficients in  line bundles $F$ that arise from various divisors $D$ on $\tilde{X}$ (i.e. $F=L_D$  for various divisors  $D$).  There exists a map 

\begin{eqnarray}\label{eq:corr}
L^{p,q}_{(2)}(U,\,\mathcal{O}(D))&\to& L^{p,q}_{(2)}(U,\,L_D)\\
f&\to& f\otimes s\notag
\end{eqnarray}

\noindent which is easily seen to be a bicontinuous isomorphism between $L^{p,q}_{(2)}(U,\,\mathcal{O}(D)))\cong L^{p,q}_{(2)}(U,\,L_D)$. The inverse to the above map is given by sending an $ L^{p,q}_{(2)}(U, L_D)\ni A \to A\cdot s^{-1}$, where $s,\,s^{-1}$ were defined in the first paragraph of section 2.2.  Based on this remark, in subsequent sections we will be tacitly identifying $H^{p,q}_{(2)}(U,\,\mathcal{O}(D))$ and $H^{p,q}_{(2)}(U,\,L_D)$. 
 
\medskip
\noindent In section 4 of the paper we shall need a generalized density lemma and closed-range property for $\overline\partial_{-D}$ (i.e. the $\overline\partial$ operator associated to the line bundle $L^*_D\cong L_{-D}$
for some divisor $D$). To simplify notation, we will consider a holomorphic line bundle $F$ over $\tilde{X}$ and a hermitian metric $h$ on it that is smooth up to $\overline{U}$. Consider the $\overline\partial_F,\,\overline\partial^*_F$ operators, defined in an 
analogous manner as before. 

\begin{lemma} The space $C^{\infty}(\overline{U},\,F)\cap \text{Dom}(\overline\partial^*_{F,h})$ is dense in the $\text{Dom}(\overline\partial_F)\cap \text{Dom}(\overline\partial^*_{F,\,h}) \cap L^{p,q}_{(2)}(U,\,F)$ for the graph norm $A\to \|A\|+\|\overline\partial_F\,A\|+\|\overline\partial^*_{F,\,h}\,A\|$. 
\end{lemma}

\medskip
\noindent {\it{Proof.}} By a partition of unity argument, it is enough to consider sections supported by $\overline{U}\cap V$, where $V$ is a small cooordinate chart over which we have a local holomorphic trivialization
$e$ of $F$. Writing $h(e,e)=e^{-\psi}$ on $V$, we see that $\overline\partial_F\,(u\otimes e)=\overline\partial u\otimes e$ and $\vartheta_{F,\,h}(u\otimes e)=(\vartheta_{\psi}u)\otimes e$, where $\vartheta_{\psi}\,u:=\vartheta u-
\partial \psi \lrcorner u$ is the formal adjoint of $\overline\partial$ with respect to the weighted $L^2$-inner product $(f,g)_{\psi}:=\int <f,g>\,e^{-\psi}\,dV$. We see that $u\otimes e\in \text{Dom}(\overline\partial^*_{F,\,h})$ if and only if $u\in \text{Dom}(\overline\partial^*)$, and then the result follows from the ordinary density lemma for scalar-valued forms.\;\;\;Q.E.D.

\medskip
\noindent Let us consider the following complex

$$ L^{p,q-1}_{(2)} (U,\,F)\overset{\overline\partial_F}\to L^{p,q}_{(2)}(U,\,F)
\overset{\overline\partial_F}\to  L^{p,q+1}_{(2)}(U,\,F).$$

\noindent  

\medskip
\noindent Recall that $U$ is a smoothly bounded strongly pseudoconvex domain in $\tilde{X}$, the Hilbert spaces are taken using the metric $h$  and  $\overline\partial^*_{F,\,h}$  denotes the Hilbert space adjoint of $\overline\partial_F$. We want to show that 

\begin{lemma} The $\text{Range}\,(\overline\partial_F)$ is closed in $L^{p,q}_{(2)}(U,\,F)$, if $q>0$.
\end{lemma}

\medskip\noindent{\it{Proof.}}  For any element $A\in \mathcal{D}_F:=\text{Dom}(\overline\partial_F)\cap \text{Dom}(\overline\partial^*_{F,\,h})\subset L^{p,q}_{(2)}(U,\,F)$ we set $|||A|||^2_{F}:=\|A\|^2+\|\overline\partial_F\,A\|^2+\|\overline\partial^*_{F,\,h} A\|^2$
where all the norms are computed with respect to $h$ and a fixed non-degenerate metric on $\tilde{X}$. The key observation in order to prove Lemma 2.2 is that if a ball in $\mathcal{D}_F$ (with respect to $|||\;|||_F$) is relatively compact in $L^{p,q}_{(2)}(U,\,F)$, then $\overline\partial_F$ has closed image in $L^{p,q}_{(2)}(U,\,F)$ and in $L^{p,q+1}_{(2)}(U,\,F)$, if $q>0$. We know that when $\partial U$ is smooth, strongly pseudoconvex and $F:=U\times \mathbb{C}$ (the scalar valued case), this observation is true (combining  Theorem 5.3.7 in \cite{ChS} and Rellich's lemma). We set $\mathcal{D}:=\text{Dom}(\overline\partial)\cap \text{Dom}(\overline\partial^*)\subset L^{p,q}_{(2)}(U)$ and $|||f|||^2
:=\|f\|^2+\|\overline\partial f\|^2+\|\overline\partial^*\,f\|^2$, in this case. Then we have the following general result: 

\begin{lemma} Let $U$ be a relatively compact subdomain  in $\tilde{X}$. 
Assume that $\{f\in \mathcal{D}:\; |||f|||\le 1\}$ is relatively compact in $L^{p,q}_{(2)}(U)$. Then $\{A\in \mathcal{D}_F:\; |||A|||_F\le 1\}$ is relatively compact in $L^{p,q}_{(2)}(U,F)$, for any holomorphic line bundle 
in a neighborhood of $\overline{U}$ and any choice of smooth metric $h$ on $F$. 
\end{lemma}

\medskip\noindent{\it{Proof.}} Cover $\overline{U}$ by relatively compact open sets $V_1,\cdots V_m$ where we have 
holomorphic trivializations $e_j$ of $F$ over $V_j$ for each $j$. Choose $\zeta_j\in C^{\infty}_0(V_j);\,0\le \zeta_j\le 1$ that form a partition of unity on $\overline{U}$. Given $A\in L^{p,q}_{(2)}(U,\,F),$ we have $s=f_j\otimes e_j$ on $U\cap V_j$ for all $1\le j\le m$. The linear map 

\begin{eqnarray*}
 \Theta: L^{p,q}_{(2)}(U,F)&\longrightarrow& \left(L^{p,q}_{(2)}(U)^{ m}\right)\\
 A&\longrightarrow& (\zeta_1 f_1^0,\zeta_2 f_2^0,\cdots,\,\zeta_m f_m^0)\\
\end{eqnarray*}

\noindent where $k^0$ denotes extension of the form $k$ by zero to $U$, is a bounded  map from $L^{p,q}_{(2)}(U,F) \to \left(L^{p,q}_{(2)}(U)^{m}\right)$, and maps $\mathcal{D}_F$ into $\mathcal{D}^{m}$. 

\medskip\noindent Let $\chi_j\in C^{\infty}_0(V_j)$ such that $\chi_j=1$ on $\text{supp}\,\zeta_j$ for all 
$j\le m$ and let us define a map $K$

\begin{eqnarray*}
 K: \left(L^{p,q}_{(2)}(U)^{m}\right)\longrightarrow L^{p,q}_{(2)}(U,F)\\
(g_1,\,\cdots g_m)&\longrightarrow& \sum_{j=1}^m (\chi_j g_j\otimes e_j)^0.\\
\end{eqnarray*}

\noindent One can easily check that $K$ is a bounded left inverse to $\Theta$. 

\medskip
\noindent Now, by elementary estimations we can show that  for all $j$ with $1\le j\le m$ we have $|||\zeta_j\,f_j|||\le C |||A|||_F$ for some 
positive constant $C$ and for all $A\in \mathcal{D}_F$. It follows that when $B$ is a $|||\;\;|||_F$-ball in $\mathcal{D}_F$, then $\Theta(B)$ is relatively compact in $L^{p,q}_{(2)}(U)^{ m}$, so $B=K(\Theta(B))$ is relatively compact in $L^{p,q}_{(2)}(U,F)$.\;\;\;Q.E.D.

\medskip\noindent One can obtain a more direct proof of Lemma 2.2 by suitably modifying H\"ormander's 
arguments in the proof of Theorem 3.4.1 in \cite{H1}. The key observation is that the assertion of the lemma is independent of a ``conformal'' change of the  metric $h$ of $F$. Setting for example $\hat{h}:=h\,\xi$ where  $\xi \in C^0(\overline{U})$ and $\xi>0$ on $\overline{U},$ would only produce equivalent norms on the Hilbert spaces that 
appear just before Lemma 2.2. Then one can use as $\xi:=e^{-\tau\,\phi}$, where $\phi$ is chosen as in the proof of 
Theorem 3.4.1 in \cite{H1} and follow H\"ormander's argument to show that the range of $\overline\partial_F$ is closed in $L^{p,q}_{(2),\,\hat{h}}(U,\,F)$ for $q>0$.

%\medskip
%\noindent 
%\medskip\noindent One can also consider $L^{p,q}_{2,\,\text{loc}}(U,\,L_D):=\{x\;\;|\;\;|x|_V<\infty\;%\text{for\;all\;open}\;V\subset\subset U\}$. If $U$ is a smoothly bounded domain with strongly pseudoconvex boundary the inclusion %$L^{p,q}_{(2)}(U,L_D)\hookrightarrow L^{p,q}_{2,\,\text{loc}}(U,L_D)$ induces isomorphisms on th ecohomology groups %$H^{p,q}_{(2)}(U,\,L_D)\cong H^{p,q}_{2,\text{loc}}(U,\,L_D)$.  

\medskip
\noindent 
\subsection{A short proof of the finite dimensionality of $L^2$-Dolbeault cohomology groups}  Let $\ell:
L^{p,q}_{2,\,\text{loc}}(\tilde{X}_r,\,\,\mathcal{O}(D))\to
L^{p,q}_{2,\,\text{loc}}(X'_r)$ be the map defined by
$\ell(g)=(\pi^{-1})^* (g)$ for $g\in
L^{p,q}_{2,\,\text{loc}}(\tilde{X}_r,\,\,\mathcal{O}(D))$. Clearly $\ell$ commutes with $\overline\partial$ and
induces a map on cohomology

$$
\ell_*: H^q(\tilde{X}_r,\,\Omega^p\otimes \mathcal{O}(D))\to
H^q(X'_r,\,\Omega^p).
$$

\medskip
\noindent 

\medskip
\noindent In section 6 of \cite{FOV} we compared various
$L^2$-$\overline\partial$-cohomology groups with certain sheaf
cohomology groups. We considered the natural inclusion $j:
L^{p,q}_{(2)}(X'_r)\to L^{p,q}_{2,\,\text{loc}}(X'_r)$ and studied
the corresponding induced homomorphism $j_*:
{H}^{p,q}_{(2)}(X'_r)\to H^q (X'_r, {\Omega^p}_{\upharpoonright
(X'_r)})$.

\begin{theorem}(Corollary 1.6 in \cite{FOV})  Let $j_*: {H}^{p,q}_{(2)}(X'_r)\to H^q
(X'_r, {\Omega^p}_{\upharpoonright (X'_r)})$ be the obvious
homomorphism induced by the inclusion $j: L^{p,q}_{(2)}(X'_r)\to
L^{p,q}_{2,\,\text{loc}}(X'_r)$.  Then the map $j_*$ is injective
for $p+q<n$ and $q>0$ and bijective for $p+q\le n-2$ and $q>0$.
\end{theorem}

\medskip
\noindent {\it{Proof of finite dimensionality of Dolbeault
cohomology groups.}} For a form $(p,q)$ form $f$ defined on $X'_r$
and square-integrable with respect to the ambient metric, its
pull-back $\pi^* f$  need not belong to
$L^{p,q}_{2,\,\sigma}(\widetilde{X_r}\setminus E)$ where $\sigma$ is
any non-degenerate metric on $\tilde{X}$. However, given $f\in
L^{p,q}_{(2)}(X'_r)\cap \text{Dom}\,(\overline\partial)$ we can
show, using  lemma 3.1 in \cite{FOV1} (comparison estimates of weighted $L^2$-norms between forms and 
their pull-backs under resolution of singularities maps), 
that $\pi^*\,f\in
\Gamma(\tilde{X}_r,\,\mathcal{L}^{p,q}(\mathcal{O}(D)))$  for some
divisor $D=\sum_{i=1}^m d_i E_i$ supported on $E$ whenever
$d_1,\cdots d_m$ are chosen large enough ($d_i>>1$ for all
$i=1,\,\cdots, m$). In addition,\; $\overline\partial \pi^*
f=\pi^*(\overline\partial f)$ for any $f\in L^{p,q}_{(2)}(X'_r)\cap
\text{Dom}\,(\overline\partial)$ and any $p,\,q\ge 0$. Taking into
account all these we  obtain a commutative diagram

$$
\xymatrix{ H^{p,q}_{(2)}(X'_r) \ar[r]^(.40){\pi^*} \ar[dr]^{j_*}
&H^q(\tilde{X}_r,\,\Omega^p\otimes
\mathcal{O}(D))\ar[d]^{\ell_*} \\
^{} &H^q(X'_r,\,\Omega^p).\\
}
$$

\noindent

\medskip
\noindent By Theorem 2.4 we know that the map $j_*$ is injective
for $p+q\le n-1,\;q>0$.  Hence the map $\pi^*$ is injective for
such $p,\,q$. As $\tilde{X}_r$ is a smoothly bounded strongly
pseudoconvex domain the cohomology groups
$H^q(\tilde{X}_r,\,\Omega^p\otimes \mathcal{O}(D))$ are finite
dimensional. Hence $H^{p,q}_{(2)}(X'_r)$ are finite dimensional
for all $p+q\le n-1,\;\,q>0$.

\medskip
\noindent Global identifications of the
$L^2$-$\overline\partial$-cohomology groups on projective surfaces
with isolated singularities with cohomology groups of appropriate
sheaves on the desingularized manifolds have been obtained by Pardon
(for cones over smooth projective curves) in \cite{P}, by Pardon and
Stern in \cite{PS1} for
$L^{n,q}_{(2)}$-$\overline\partial$-cohomology groups of projective
varieties with arbitrary singularities and recently by Ruppenthal
\cite{Rup3} for a large class of compact pure dimensional Hermitian
complex spaces with isolated singularities.

\section{Proof of Theorem 1.1}

\medskip
\noindent \subsection{Exceptional sets.} We  shall recall the notion
of exceptional sets (in the sense of Grauert \cite{Gr}) and some key
results regarding these sets that will be needed in the paper.

\medskip
\noindent{\bf{Definition.}} Let $X$ be a complex space. A compact
nowhere discrete, nowhere dense analytic set $A\subset X$ is
exceptional if there exists a proper, surjective map $\pi:X\to Y$
such that $\pi(A)$ is discrete, $\pi: X\setminus A \to Y\setminus
\pi(A)$ is biholomorphic and for every open set $D\subset Y$ the map
$\pi^*: \Gamma(D, \mathcal{O})\to \Gamma
(\pi^{-1}(D),\,\mathcal{O})$ is surjective.

\smallskip
\noindent We usually say that $\pi$ collapses or blows down $A$.

\smallskip
\noindent If $V$ is a Stein neighborhood of $\pi(A)$ then
$\pi^{-1}(V)$ is a $1$-convex space with maximal compact analytic
set $A$ and $\pi_{\upharpoonright \pi^{-1}(V)}$ is the Remmert
reduction.

\medskip
\noindent {\bf{Characterization of exceptional sets}}

\medskip
\noindent Below we collect some basic results regarding exceptional
sets.

\medskip
\noindent $\alpha)$ (Theorem 4.8, page 57 in \cite{Lauf2}) Let $X$
be an analytic space and $A$ a compact, nowhere discrete analytic
subset. $A$ is exceptional if and only if there exists a
neighborhood $U$ of $A$ such that the closure of $U$ in $X$ is
compact, $U$ is strictly Levi pseudoconvex and $A$ is the maximal
compact analytic subset of $U$. Also, $A$ is exceptional if and only
if $A$ has arbitrarily small strictly pseudoconvex neighborhoods.

\medskip
\noindent $\beta)$ (Lemma 3.1 in \cite{Lauf1}) Let $\pi: U\to Y$
exhibit $A$ as exceptional in $U$ with $Y$ a Stein space. If
$U\supset V$ with $V$ holomorphically convex neighborhood of $A$ and
$\mathcal{F}$ is a coherent analytic sheaf on $U$, then the
restriction map $\rho: H^i(U,\,\mathcal{F})\to H^i(V,\,\mathcal{F})$
is an isomorphism for $i\ge 1$.

\subsection{Local cohomology along exceptional sets} In this section
we recall Stephen Yau's and Karras' results  that describe the local
cohomology along exceptional sets. Our earlier work on Hartogs'
extension theorems on Stein spaces (see \cite{OV2}) indicated to us
the importance of the local cohomology exact sequences and led us to
the discovery of these theorems.

\medskip
\noindent For a sheaf of abelian groups $\mathcal{F}$ on a
paracompact, Hausdorff space $X$ and for $K$ a closed subset of $X$,
let $\Gamma_K (X,\mathcal{F})$ denote the sections on $X$ with
support in $K$. Consider a flabby resolution of $\mathcal{F};\;\;\;
0\to \mathcal{F}\to \mathcal{C}^{0}\overset{d_0}\to
\mathcal{C}^{1}\overset{d_1}\to \cdots$. The cohomology groups with
support in $K$ are defined by $H^{*}_K (X, \mathcal{F}):=H^{*}
(\Gamma_K \left(X, \mathcal{C}^{\bullet})\right)$, i.e. they are the
cohomology groups of the complex $\left(\Gamma_K(X,
\mathcal{C}^{k}), d_k \right)$. Since each $\mathcal{C}^{\bullet}$
is flabby, we have a short exact sequence $0\to \Gamma_K (X,
\mathcal{C}^{\bullet})\to \Gamma (X, \mathcal{C}^{\bullet})\to
\Gamma (X\setminus K, \mathcal{C}^{\bullet})\to 0$. This induces a
long exact sequence on cohomology

\begin{equation}\label{eq:oles}
0\to H^0_K(X,\mathcal{F})\to H^0(X, \mathcal{F})\to H^0(X\setminus
K, \mathcal{F})\to H^1_K(X,\mathcal{F})\to H^1(X,\mathcal{F})\to...
\end{equation}

\noindent It is a standard fact from sheaf cohomology theory that
$H^i_K(X,\mathcal{F})\cong H^i_K(U, \mathcal{F})$ where $U$ is an
open neighborhood of $K$ in $X$. The fact that
$H^i_K(U,\mathcal{F})$ is independent of the neighborhood $U$ of $K$
is referred to as excision.

\medskip
\noindent On the other hand, one could also consider the cohomology
with compact support on $X$ and define $\Gamma_c(X,\mathcal{F})$ to
be the group of global sections of $\mathcal{F}$ whose supports are
compact subsets of $X$. Let $Y$ be a compact subset of $X$. Letting
$\{\mathcal{C}^i\}$ denote the canonical resolution of
$\mathcal{F}$, we have an inclusion of complexes

$$
\left(\Gamma_Y(X,\,\mathcal{C}^{\cdot}),\,d_{\cdot}\right)\hookrightarrow
\left(\Gamma_c(X,\,\mathcal{C}^{\cdot}),\,d_{\cdot}\right)
$$

\noindent which induces natural homomorphisms $\gamma_i:
H^i_Y(X,\,\mathcal{F})\to H^i_{c}(X,\,\mathcal{F})$ for all $i\ge
0$.

\medskip
\noindent In general we do not have enough information on the maps
$\gamma_i$ but in the special case of exceptional sets we can obtain
very precise information about them. For the remainder of this
section let $X$ be a reduced complex space and $Y=E$ be an
exceptional subset of $X$. Then we know that there exists a strongly
pseudoconvex neighborhood $M$ of $E$ in $X$ and a non-negative
exhaustion function $\phi$ on $M$ such that $\phi$ is strongly
plurisubharmonic on $M\setminus E$ and $E=\{x\in M;| \;\phi(x)=0\}$.
By excision, $H^i_E(X,\mathcal{F})=H^i_E (M,\,\mathcal{F})$ for all
$i\ge 0$ and therefore we have natural homomorphisms $\gamma_i:
H^i_E (X,\,\mathcal{F})\to H^i_c(M,\,\mathcal{F})$ for all $i\ge 0$.
Karras showed that under circumstances these maps $\gamma_i$ are
isomorphisms for some $i$.

\begin{theorem} (Proposition 2.3 in \cite{Kar}) Let $X$ be a reduced
complex space and $E$ an exceptional subset of $X$. If $\mathcal{F}$
is a coherent analytic sheaf on $X$ such that
$\text{depth}_x\,\mathcal{F}\ge d$ for $x\in M\setminus E$, then

$$
\gamma_i: H^i_E(X,\,\mathcal{F})\to H^i_c(M,\,\mathcal{F})
$$

\noindent is an isomorphism for $i<d$.
\end{theorem}

\medskip
\noindent

\noindent Once we have Theorem 3.1, we can very easily obtain the
following corollary:

\begin{corollary} Let $E$ be an
exceptional set of an $n$-dimensional complex manifold $M$. Then

$$
H^i_E(M,\,\mathcal{O}_M)=0\phantom{asdasda}\text{for}\phantom{asdjhja}i<n.
$$
\end{corollary}

\medskip
\noindent{\it{Proof.}} Since $E$ is an exceptional set of $M$, $M$
is a strongly pseudoconvex manifold and let $p:M\to S$ denote the
Remmert reduction map. For every coherent analytic sheaf
$\mathcal{F},\;\; H^{i}(M,\,\mathcal{F})$ are finite dimensional
for $i>0$. Hence we can apply Serre's duality theorem for
$\mathcal{F}=\omega_M=\Omega^n$, the sheaf of holomorphic
$n$-forms on $M$. Then $H^{n-i}(M,\omega_M)\cong
H^{i}_c(M,\,\mathcal{O})$ for all $i<n$. Since $R^i p_*\omega_M=0$
for all $i>0,$ by Takegoshi's relative vanishing theorem in
\cite{Tak}, we have
$H^{n-i}(M,\omega_M)=H^{n-i}(S,\,p_*(\omega_M))$. But the latter
cohomology groups vanish since $S$ is Stein, $p_*(\omega_M)$ is
coherent and $n-i>0$. Therefore $H^i_c(M,\,\mathcal{O}_M)=0$ for
all $i<n$. Since $M$ is a manifold
$\text{depth}_x(\mathcal{O}_M)=n$ for all $x\in M\setminus E$;
hence we can apply Theorem 3.1 to conclude that
$H^i_E(M,\,\mathcal{O}_M)\cong H^i_c(M,\,\mathcal{O}_M)=0$ for all
$i<n$.

\medskip
\noindent
\subsection{Proof of Theorem 1.1} Suppose now that $X$ is a pure  $n$-dimensional
($n\ge 2$) complex analytic set in $\mathbb{C}^N$ with an isolated
singularity at $0$ and let $X_r=X\cap B_r$ be a small Stein
neighborhood of $0$ with smooth boundary. Let $\pi: \tilde{X}\to X$
be a desingularization of $X$.  Then  $E:=\pi^{-1}(0)$ (the
exceptional locus of the desingularization) is an exceptional set in
the sense of Grauert and let $\tilde{X}_r:=\pi^{-1}(X_r)$.  Let
$\sigma$ be a positive definite metric on $\tilde{X}$. In what
follows $\mathcal{L}^{\bullet,\,\bullet}$ represents the sheaves of
differential forms that were introduced in section 2.1.

\medskip
\noindent Let $r>0$ be a regular value of $\|\;\;\|\circ \pi$ on
$\tilde{X}_R$ with $0<r<R$, $\tilde{X}_r$ is a relatively compact
domain with smooth strongly pseudoconvex boundary in $\tilde{X}_R$.
It is a standard fact that the inclusion of the following complexes

$$
L^{0,^{.}}_{(2)}(\tilde{X}_r)\cap
\mathcal{D}(\overline\partial)\overset{\tilde{j}}\hookrightarrow
\mathcal{L}^{0,\,.}(\tilde{X}_r)
$$

\noindent induces isomorphisms on the corresponding cohomology
groups $H^{0,q}_{(2)}(\tilde{X}_r)\cong
H^q(\tilde{X}_r,\,\mathcal{O})$ for $q>0$.

\medskip
\noindent  By Theorem 2.4, we know that for $1\le q\le n-2$ we
have ${H}^{0,q}_{(2)}(X'_r)\cong H^q(X'_r, \mathcal{O})$. The
latter sheaf cohomology groups are isomorphic to
$H^q(\tilde{X}_r\setminus E,\,\mathcal{O}_{\tilde{X}_r})$.
Consider the long exact local cohomology sequence

\begin{eqnarray}\label{eq:2lelcs}
....\to H^i_E(\tilde{X}_r,\,\mathcal{O})\to
H^i(\tilde{X}_r,\,\mathcal{O})\overset{r_*} \to
H^i(\tilde{X}_r\setminus E,\,\mathcal{O})\to
H^{i+1}_E(\tilde{X}_r,\,\mathcal{O})\to...
\end{eqnarray}

\medskip
\noindent If $1\le q\le n-2$, then by Corollary 3.2 we have
$H^q_E(\tilde{X}_r,
\,\mathcal{O})=H^{q+1}_E(\tilde{X}_r,\,\mathcal{O})=0$; hence from
$(\ref{eq:2lelcs})$ we can conclude that
$H^q(\tilde{X}_r,\mathcal{O})\cong H^q(\tilde{X}_r\setminus
E,\,\mathcal{O})\cong H^q(X'_r,\,\mathcal{O})\cong
H^{0,q}_{(2)}(X'_r)$.

\medskip
\noindent We shall construct now the map $\phi_*$ that appears in
Theorem 1.1. Recall that for a $(0,q)$ form in $\tilde{X}_r$, we
have $\|g\|_{L^2_{\gamma}(\tilde{X}_r)}\le
C\,\|g\|_{L^2(\tilde{X}_r)}$, where $\gamma$ is the
``pseudometric'' that arises from the pull-back of the Euclidean
metric in $X'$, since $|\,\;|_{\gamma}\le C_0\,|\,\;|_{\sigma}$
and $p=0$. Moreover
$\|(\pi^{-1})^*\,g\|_{L^2(X'_r)}=\|g\|_{L^2_{\gamma}(\tilde{X}_r)}$.
Now, for $g\in L^{0,q}_{(2)}(\tilde{X}_r)$ we have
$\|(\pi^{-1})^*\,g\|_{L^2(X^*_r)}\le C\,\|g\|_{L^2(\tilde{X}_r)}$,
since we pass to a smaller norm. Hence for all $q$ with $0\le q\le
n$, there exists a bounded linear map:

\begin{eqnarray*}\label{eq:dfphi}
\phi: L^{0,q}_{(2)}(\tilde{X}_r) &\to &L^{0,q}_{(2)}(X'_r)\\
u&\to& (\pi^{-1})^* (u).\\
\end{eqnarray*}

\medskip
\noindent Then we have a commutative diagram of complexes

$$
\xymatrix{ \mathcal{L}^{0,\,.}(\tilde{X}_r)\ar[drr]^{r}  \\
&L^{0,\,\cdot}_{(2)}(\tilde{X}_r)\cap \mathcal{D}(\overline\partial)
\ar[ul]\ar[d]^{\phi}\ar[r] &   \mathcal{L}^{0,\,\cdot}(\tilde{X}_r\setminus E)\ar[d]^{\cong}\\
 &L^{0,\,\cdot}_{(2)}(X'_r)\cap \mathcal{D}(\overline\partial)
 \ar[r]^{j} & \mathcal{L}^{0,\,.}(X'_r)\\
 }
$$

\noindent which induces the following commutative diagram:
\medskip
\noindent
$$
\xymatrix{H^q(\tilde{X}_r,\,\mathcal{O})\ar[drr]^{r_*}
\\
&H^{0,q}_{(2)}(\tilde{X}_r) \ar[r] \ar[ul]^{\cong} \ar[d]^{\phi_*}
& H^q(\tilde{X}_r\setminus E,\,\mathcal{O}) \ar[d]^{\cong}\\
&  H^{0,q}_{(2)}(X'_r) \ar[r]^{j_*}    & H^q(X'_r,\,\mathcal{O}). }
$$

\medskip
\noindent Since for $1\le q\le n-2$\;\; $r_*$ and $j_*$ are
isomorphisms, the commutativity of the above diagram will imply that
$\phi_*$ is an isomorphism for $1\le q\le n-2$. On the other hand
for $q=n-1$\; the maps $r_*$ and $j_*$ are only injective, hence
$\phi_*$ is an injective map.

\medskip
\noindent 

\section{Proof of Theorem 1.2}

\subsection{Different extensions for
$\overline\partial$-operator} For the proof of Theorem 1.2 another
closed extension of the $\overline\partial$-operator will play a
key role. Let $\overline\partial^1$ denote the graph closure in
$L^2$ of $\overline\partial$ acting on forms $f$ with coefficients in
$C^{\infty}_{0}(\overline{X}_r\setminus \{0\})$. We can also
consider the minimal extension of the $\overline\partial$-operator
on $X'$. More precisely we let $\overline\partial_{\text{min}}$
denote the graph closure in $L^2$ of $\overline\partial$ acting on
forms with coefficients in $ C^{\infty}_{0}({X}_r\setminus \{0\})$
(Dirichlet conditions on both the boundary of $X_r$ and the
singularity $0$). It is easy to check that

\begin{lemma}
$\text{Dom}(\overline\partial^1)=\{f\in
\text{Dom}(\overline\partial):\;\chi\,f\in
\text{Dom}(\overline\partial_{min})\;\;\;\text{for\;\; a\;\;
cut-off\;\; function}\; \chi\in C^{\infty}_0(X_r)\;\}$.
\end{lemma}

\medskip
\noindent Forms of bidegree $(0,n-1)$ in
$\text{Dom}(\overline\partial)(\tilde{X}_r)$ are in the domain of
$\overline\partial^1$. More precisely, we have:

\begin{lemma} If $h\in
\text{Dom}(\overline\partial)\cap L^{0,n-1}_{(2)}(\tilde{X}_r),$
then $\phi(h)\in \text{Dom}(\overline\partial^1)$.
\end{lemma}

\noindent {\it{Proof.}} We will distinguish two cases:

\medskip
\noindent {\it{Case I.}} Let $h\in
\text{Dom}(\overline\partial_{\tilde{X}_r})\cap
L^{0,n-1}_{(2)}(\tilde{X}_r)$ and smooth in $\tilde{X}_r$ (thus
bounded near $E$). By a partition of unity argument we can assume
that the support of $h$ is contained in a coordinate domain $U$,
where $U\overset{\cong}\to U'\subset \mathbb{C}^n$ when $E\cap
U=\{z\in U: \;z_1\,z_2\cdots z_m=0\}$ for some $m$ with $1\le m\le
n$. We choose a family of cut-off functions $\chi_k$ that satisfy:
i) $\chi_k(z)=1$ when $\text{dist}(z,\,E)\ge \frac{1}{k}$ and
$\chi_k(z)=0$ near $E$, and ii) $|\overline\partial \chi_k(z)|\le
C\,k$ for all $k$. Now $\phi(h)$ has compact support $\pi(U)\cap
X_r$ and in order to show that it belongs to
$\text{Dom}(\overline\partial^{1})$ it suffices by Lemma 4.1 to
show that $\phi(h)\in
\text{Dom}\,(\overline\partial_{\text{min}})$. Since
$\overline\partial_{\text{min}}={{\vartheta^*}_{\text{max}}},$
(the Hilbert space adjoint of $\vartheta_{\text{max}}$) we must
show that

$$(\overline\partial \phi(h),\;w)=(\phi(h),\;\vartheta_{max} w)$$

\noindent for all $w\in \text{Dom}(\vartheta_{max})\cap
L^{0,n}_{(2)}(X'_r)$. Let us set $\psi_k:=\chi_k\circ \pi^{-1}$.
Then we have

$$
(\phi(h),\,\vartheta_{max}
w)=\text{lim}(\psi_k\,\phi(h),\,\vartheta_{max}
w)=\text{lim}(\overline\partial \psi_k \wedge
\phi(h),\,w)+(\overline\partial \phi(h),\,w).
$$

\noindent But

$$
|\int_{\tilde{X}_r} \overline\partial \chi_k \wedge h\wedge
\pi^*(\overline {*}\, w) \;|\le \|h\|_{\infty}
\,\left(\int_{\tilde{X}_r} |\overline\partial \chi_k|^2\,dV
\right)^{\frac{1}{2}}\;\;\left(\int_{\text{supp}\,\overline\partial\chi_k}
|\pi^* w|^2\,dV \right)^{\frac{1}{2}}=A\,\,B\,\,C
$$

\noindent where $C\to 0$ as $k\to \infty$, while $B$ is easily
seen to be uniformly bounded. Hence, $\phi(h)\in
\text{Dom}\,(\overline\partial^1)$.

\medskip
\noindent{\it{Case II.}} Let $h\in \text{Dom}(\overline\partial)\cap
L^{0,n-1}_{(2)}(\tilde{X}_r)$. Since the smooth forms in
$\overline{\tilde{X}_r}$ are dense in
$\text{Dom}(\overline\partial)$ in the graph norm, there exist
$h_{\nu}\in C^{\infty}(\overline{\tilde{X}_r})$ such that $h_{\nu}
\to h$ in the graph norm as $\nu \to \infty$. But then
$\phi(h_{\nu})\in \text{Dom}(\overline\partial^1)$ by Case I and
converge to $\phi(h)$ in the graph norm in $L^{0,n-1}_{(2)}(X'_r)$.
Recall that $\overline\partial^1$ is a closed operator, hence
$\phi(h)\in \text{Dom}(\overline\partial^1)$.

\medskip
\noindent
\subsection{Preliminaries from \cite{FOV}.}  In a previous work (lemma 3.4 in \cite{FOV}), we showed that for
$f\in Z^{p,q}_{(2)}:=L^{p,q}_{(2)}(X'_r)\cap
\text{kern}(\overline\partial)$ with $p+q\le n-1$ and $q>0$, the
equation $\overline\partial u=f$ is solvable in
$L^{p,q-1}_{(2)}(X'_r)$ if and only if the equation
$\overline\partial v=f$ is solvable in $L^{p,q-1}_{(2)}(X_r\setminus
\overline{B}_{r_0})$ for some $0<r_0<r$. In addition we showed
(Proposition 3.5 in \cite{FOV}) that the equation $\overline\partial
u=f$ is solvable in $L^{p,q-1}_{(2)}(X_r\setminus
\overline{B}_{r_0})$  with $0<r_0<r$, for $f$ in a closed subspace
of finite codimension in $Z^{p,q}_{(2)}$ when $p+q<n,\,q>0$. Let us
recall Case II in the proof of Proposition 3.5 in \cite{FOV}. Let
$f\in Z^{0,n-1}_{(2)}(X'_r)$ and let $\chi\in C^{\infty}_0(X_r)$
with $\chi=1$ near $0$ and $\text{supp}\,\chi \subset X_{\rho}$ with $0<\rho<r$. It was shown that $\overline\partial
w=\pi^*(\overline\partial\chi \wedge f)$ had a solution in
$L^{0,n-1}_{(2)}(\tilde{X}_r)$, compactly supported in $\tilde{X}_r$
if and only if

\begin{eqnarray}\label{eq:printest}
\int_{X'_{\rho}} f\wedge \overline\partial \chi \wedge \psi=0\;\;
\text{for\;\;all}\;\;\psi\in L^{n,0}_{(2)}(X'_{\rho})\cap
\text{kern}(\overline\partial):=H^{n,0}_{(2)}(X'_{\rho}).
\end{eqnarray}

\medskip
\noindent Condition $(\ref{eq:printest})$ can be derived from the following 
weaker consition:

\begin{eqnarray}\label{eq:intest}
\int_{X'_{r}} f\wedge \overline\partial \chi \wedge \psi=0\;\;
\text{for\;\;all}\;\;\psi\in L^{n,0}_{(2)}(X'_{r})\cap
\text{kern}(\overline\partial):=H^{n,0}_{(2)}(X'_r).
\end{eqnarray}

\noindent This is a consequence of the following fact:

\begin{lemma} The pair $(X'_r,\,X'_{\rho})$ is an $L^2$-Runge pair, i.e. the restriction map $r:\mathcal{O}L^2(X'_r)\to \mathcal{O}L^2(X'_{\rho})$ has dense image. 
\end{lemma}

\smallskip
\noindent{\it{Proof.}} Let $h\in \mathcal{O}L^2(X'_{\rho})$. We need to show that there exists a sequence of 
functions $h_{\nu}\in \mathcal{O}L^2(X'_r)$ such that $\|h_{\nu}-h\|_{L^2(K)}<\epsilon$, where $\epsilon>0$ and $K$ a compact subset of $X'_{\rho}$. Suppose that $K\subset A_0:=X\cap \{r_0^*<\|z\|<r_1^*\}\subset\subset A:=X\cap \{r_0<\|z\|<r_1\}$ where $0<r_0<r_0^*<r_1^*<r_1<\rho<r$. Let $\phi\in C^{\infty}_0(A)$ with $\phi=1$ on $\overline{A_0}$ and let us look at the $(0,1)$-form $g:=\overline\partial(\phi\,h)$; we can write $g=g'+g''$ where $g'$ is supported on $X\cap \{\|z\|>r_1^*\}$, $g''$ is supported in $X_{r_0^*}$ and both are 
$\overline\partial$-closed on $X'_r$. Using Proposition 3.1 from \cite{FOV}, we know that there exists a 
solution $u''$ satisfying $\overline\partial u''=g''$ on $X'$, compactly supported in $X_{r_0^*}$ and in $L^{0,0}_{(2)}(X')$. We consider a convex, increasing function $\xi\in C^{\infty}(\mathbb{R})$ with $\xi(t)=0$ if $t\le r_1^*$ and $\xi(t)>0$ if $t>r_1^*$. Let $\psi(z)=\xi(\|z\|)$. By our choice of $\xi$, we know that the $\text{min}\,\{\psi(z);\;z\in \text{supp}\,g'\}=c>0$. Applying Theorem 1.3 in \cite{FOV}, we obtain a solution $u'_{\nu}$ satisfying $\overline\partial u'_{\nu}=g'$ in $X'_r$ and 

$$
\int_{X'_r} |u'_{\nu}|^2\, e^{-\nu\,\psi}\, dV\le Ce^{-\nu\,c}\,\|g'\|^2
$$

\noindent for $\nu\ge \nu_0$. Hence, $\int_{A_0} |u'_{\nu}|^2\,dV\le C\, e^{-c\,\nu}\,\|g'\|^2$. 

\smallskip\noindent Let $h_{\nu}:=\phi\,h-u'_{\nu}-u''$. Then $h_{\nu}\in L^{0,0}_{(2)}(X'_r), \,
\overline\partial h_{\nu}=0$ on $X'_r$ and 

$$
\int_K |h_{\nu}-h|^2\, dV<\epsilon
$$

\noindent for $\nu$ sufficiently large. Q.E.D. 

\medskip
\noindent

\medskip
\noindent Condition $(\ref{eq:intest})$ is independent of the
choice of the cut-off function $\chi$. Also, if
$f=\overline\partial u$ near the support of
$\overline\partial\chi$ then $f$ satisfies $(\ref{eq:intest})$,
since

$$
\int_{X'_{r}} f \wedge \overline\partial\chi \wedge \psi=\int_{X'_{r}}
d(u\wedge \overline\partial\chi\wedge \psi)=0$$

\noindent by Stokes' theorem. Hence condition $(\ref{eq:intest})$
depends only on the equivalence class $[f]\in
H^{0,n-1}_{(2)}(X'_r)$.

\medskip
\noindent Let $\mathcal{M}:=\{f\in
Z^{0,n-1}_{(2)}(X'_r);\;\;f\;\;\text{satisfies}\;(\ref{eq:intest})\;\}.$

\medskip
\noindent If $f\in \mathcal{M}$, we can write $f=\phi
(w+\pi^*((1-\chi)\,f)))+(\chi f-\phi(w))$ -where  $w$ is the square-integrable,  compactly supported form in $\tilde{X}_r$ that satisfies $\overline\partial w=\pi^*(\overline\partial\,\chi\wedge f)$.  Each
term to the right-hand side of the previous equation is
$\overline\partial$-closed and the second one has compact support
in $X_r$, hence it is $\overline\partial$-exact by Proposition 3.1
in \cite{FOV} (which is an $L^2$-solvability result for
square-integrable, $\overline\partial$-closed forms with compact support in $X_r$). 
Therefore we can write $[f]=\phi_*([w+\pi^*((1-\chi)\,f)])$.

\medskip
\noindent On the other hand, if $f=\phi(g)$ for some $g\in
L^{0,n-1}_{(2)}(\tilde{X}_r)\cap \text{kern}(\overline\partial),$
Lemma 4.2 from section 4.1 tell us that $\phi(g)$ will belong in the
domain of $\overline\partial^1$ and $\overline\partial^1 \phi(g)=0$.
Hence there exist $h_{\nu}\in C^{\infty}_{0}(\overline{X}_r\setminus
\{0\})$ such that $h_{\nu}\to \phi(g)$ and $\overline\partial
h_{\nu}\to  0$ in $L^2$. The latter would imply for $\psi\in
H^{n,0}_{(2)}(X'_{r})$ that

\begin{eqnarray}\label{eq:d1van}
\int_{X'_{r}} f \wedge \overline\partial\chi \wedge \psi=\lim_{\nu\to
\infty} \int_{X'_{r}} h_{\nu}\wedge \overline\partial \chi \wedge
\psi= (-1)^{n-1}\,\lim_{\nu \to \infty} \int_{X'_{r}} d(h_{\nu} \wedge
\chi\,\wedge \psi)=0,
\end{eqnarray}

\noindent by Stokes' theorem on $X'_{r}$.

\medskip
\noindent We have thus shown the following

\begin{lemma} $\mathcal{M}=
\phi(Z^{0,n-1}_{(2)}(\tilde{X}_r))+\text{Im}\,\overline\partial^{0,n-2}(X'_r)$.
\end{lemma}

\medskip
\noindent In \cite{FOV}, we showed that $\int_{X'_{\rho}}
\overline\partial \chi \wedge f\wedge \psi=0$ is satisfied for all
$f\in Z^{0,n-1}_{(2)}(X'_r)$ when $\psi\in \Gamma (X_{\rho},
\overset{\circ} \omega)$ for some coherent $\mathcal{O}_X$-module
$\overset{\circ}\omega$. The module $\Gamma(X_{\rho},\,\overset{\circ}
\omega)$ was shown to have finite codimension in
$\Gamma(X_{\rho},\,\omega)$-where $\omega$ was Grothendieck's dualizing
sheaf-and a fortiori in $L^{0,\,n}_{2,\,\text{loc}}(X_{\rho})\cap
\text{kern}(\overline\partial)$. Now, if $a_1,\,\cdots a_m$ span the
complementary subspace to $\Gamma(X_{\rho},\,\overset{\circ}\omega)$ we
see that

$$
\mathcal{M}=\{f\in Z^{0,n-1}_{(2)}(X'_r);\; \int_{X'_{\rho}} f \wedge
\overline\partial \chi\wedge
a_j=0,\;\;\;\text{for\;all}\;j=1,\cdots,m\}.
$$

\noindent Hence the codimension of $\mathcal{M}$ in
$Z^{0,n-1}_{(2)}(X'_r)$ is at most $m$. In what follows we will
identify the subspace of $H^{n,0}_{(2)}(X'_{r})$ for which
$\int_{X'_{r}} f \wedge \overline\partial \chi \wedge \psi=0$ for
all $f\in Z^{0,n-1}_{(2)}(X'_r)$.

\medskip
\noindent
\subsection{Construction of a non-degenerate pairing.} The
construction in section 4.2 allow us to consider the following
pairing: Let $\chi\in C^{\infty}_0(X_r)$ such that $\chi=1$ in a
neighborhood of $0$ and $\text{supp}\chi\subset X_{\rho}$ with $0<\rho<r$. 
Take a pair $(f,\,\psi)\in
Z^{0,n-1}_{(2)}(X'_r)\times H^{n,0}_{(2)}(X'_{r})$ and assign to it
the number

\begin{eqnarray}\label{eq:ndp}
<f,\,\psi>=\int_{X'_{r}} f \wedge \overline\partial \chi \wedge \psi.
\end{eqnarray}

\medskip
\noindent Recall from our discussion above that $<f,\psi>=0$ for all
$\psi\in H^{n,0}_{(2)}(X'_{r})$ is equivalent to the fact that $f\in
\mathcal{M}$ which in its turn is equivalent to the fact that
$f=\phi({g})+\overline\partial u$ where ${g}\in
L^{0,n-1}_{(2)}(\tilde{X}_r)\cap \text{kern}(\overline\partial)$
(hence $\phi(g)\in \text{Dom}\,(\overline\partial^1)$) and $u\in
L^{0,n-2}_{(2)}(X'_r)$, compactly supported in $X_r$.

\medskip
\noindent
\begin{proposition} We have $<f,\,\psi>=0$ for all $f\in
Z^{0,n-1}_{(2)}(X'_r)$ if and only if $\psi\in
\text{kern}\,({\overline\partial}^1)_{X'_r}$.
\end{proposition}

\medskip
\noindent {\bf{Remark:}} The first paragraph in section 4.3 and
Proposition 4.4 will allow us to say that $(\ref{eq:ndp})$ is a
non-degenerate pairing from

$$
\dfrac{Z^{0,n-1}_{(2)}(X'_r)}{\phi(Z^{0,n-1}_{(2)}(\tilde{X}_r))+\text{Im}\,\overline\partial^{0,n-2}}
\times
\dfrac{\text{kern}{\overline\partial}^{n,0}_{X'_{r}}}{\text{kern}{\overline\partial^1}^{n,0}_{X'_{r}}}
\to
\mathbb{C}
$$

\noindent or equivalently

$$
\dfrac{H^{0,n-1}_{(2)}(X'_r)}{\phi_*(H^{0,n-1}(\tilde{X}_r))} \times
\dfrac{\text{kern}(\overline\partial)^{n,0}}{\text{kern}(\overline\partial^1)^{n,0}}\to
\mathbb{C}
$$

\medskip
\noindent {\bf{Remark:}} Due to the injectivity of $\phi_*$ one can obtain the following bound on 
the complex dimension of $H^{0,n-1}_{(2)}(X'_r)$

\begin{eqnarray}\label{eq:dimform}
\text{dim}_{\mathbb{C}}
H^{0,n-1}_{(2)}(X'_r)=\text{dim}_{\mathbb{C}}
H^{0,n-1}(\tilde{X}_r)+\text{dim}_{\mathbb{C}}
\dfrac{\text{kern}(\overline\partial)^{n,0}}{\text{kern}(\overline\partial^1)^{n,0}}.
\end{eqnarray}

\medskip
\noindent{\it{Proof of Proposition 4.5.}} $\Leftarrow$ If $\psi\in
\text{kern}({\overline\partial^1})_{X'_{r}}$, then there exist $\psi_{\nu}\in
C^{\infty}(\overline{X}_r\setminus \{0\})$ such that $\psi_{\nu}\to
\psi$ in $L^2$ and $\overline\partial \psi_{\nu}\to 0$ in $L^2$ as
$\nu\to \infty$. But then $\int_{X'_{r}} f\wedge \overline\partial
\chi \wedge \psi=\lim_{\nu\to \infty} \int_{X'_{r}} f\wedge
\overline\partial \chi \wedge \psi_{\nu}=\lim_{\nu\to \infty}\,
[\,\int_{X'_{r}} f\wedge \overline\partial(\chi
\,\psi_{\nu})-\int_{X'_{r}} f\chi\wedge \overline\partial
\psi_{\nu}\,]=(-1)^{n} \lim_{\nu\to \infty} \int_{X'_{r}}
\overline\partial f \wedge \chi\,\psi_{\nu}-0=0.$

\medskip
\noindent $\Rightarrow$ Let us assume that $<f,\,\psi>=0$ for all
$f\in Z^{0,n-1}_{(2)}(X'_r)$. We want to show that $\psi\in
\text{kern}({\overline\partial^1})_{X'_{r}}$. It suffices to show that $\psi\in
\text{Dom}(\overline\partial^1)_{X'_{r}}$. By Lemma 4.1, this is equivalent
to showing that $\chi\, \psi \in
\text{Dom}(\overline\partial_{\text{min}})$. Recall that
$(\overline\partial_{\text{min}})^*=\vartheta_{\text{max}}$. Hence
to show that $\chi\, \psi \in
\text{Dom}(\overline\partial_{\text{min}})$ it would suffice to show
that $\chi\,\psi\in\text{Dom}(\vartheta_{\text{max}})^*$  or
equivalently

\begin{eqnarray}\label{eq:L2cond}
(\overline\partial
(\chi\psi),\,g)=(\chi\psi,\,\vartheta_{\text{max}}g)
\end{eqnarray}

\noindent for all $g\in
\text{Dom}(\vartheta_{\text{max}})^{n,1}=\{g\in
L^{n,1}_{(2)}(X'_r);\;\;\vartheta_{\text{max}}\,g\in
L^2\;\;(\text{weakly})\;\}$.

\medskip
\noindent The operator  $\overline {*}: L^{n,1}_{(2)}(X'_r)\to
L^{0,n-1}_{(2)}(X'_r)$ is an isometry mapping from
$\text{Dom}(\vartheta_{\text{max}}) \to
\text{Dom}(\overline\partial)$ (here $\overline\partial$ denotes the
maximal (weak) extension). Hence $(\ref{eq:L2cond})$ is equivalent
to

\begin{eqnarray}\label{eq:newL2}
\int_{X'_r} \overline\partial (\chi\,\psi)\,\wedge
w=-(-1)^n\,\int_{X_r}  \chi\,\psi\,\wedge \overline\partial w
\end{eqnarray}

\noindent for all $w\in \text{Dom}(\overline\partial)\cap
L^{0,n-1}_{(2)}(X'_r)$.

\medskip
\noindent Clearly $(\ref{eq:newL2})$ holds for all $w\in
Z^{0,n-1}_{(2)}(X'_r)$, by the assumption $<f,\,\psi>=0$ for all
$f\in Z^{0,n-1}_{(2)}(X'_r)$. Let us consider an arbitrary element
$w \in  \text{Dom}(\overline\partial)\cap L^{0,n-1}_{(2)}(X'_r)$.
Then $\overline\partial w \in Z^{0,n}_{(2)}(X'_r)$ and
$\pi^*(\overline\partial w)\in H^{0,n}_{(2)}(\tilde{X}_r)$ (a few
words are in order here: a)\;$\pi^*(\overline\partial w)\in
L^{0,n}_{(2),\,\gamma}(\tilde{X_r})=L^{0,n}_{(2),\,\sigma}$ where
$\sigma$ is a non-degenerate metric on $\tilde{X}_r$ and $\gamma$ is
the pull-back of the Euclidean metric under $\pi$. b)
$\pi^*(\overline\partial w)$ is $\overline\partial$-closed in
$\tilde{X}_r\setminus E$, but since $\pi^*(\overline\partial w)\in
L^{0,n}_{(2),\,\sigma}(\tilde{X_r})$ it can be extended as a
$\overline\partial$-closed form in $\tilde{X}_r$. We shall still
denote the extended form as $\pi^*(\overline\partial w)$.) Now, as
$\tilde{X}_r$ is a smoothly bounded domain with strongly
pseudoconvex boundary we have $H^{0,n}_{(2)}(\tilde{X}_r)\cong
H^n(\tilde{X}_r,\,\mathcal{O})=0$ (the latter due to work of Siu,
\cite{Siu}). Hence, there exists a solution $g\in
L^{0,n-1}_{(2)}(\tilde{X}_r)$ such that $\overline\partial
g=\pi^*(\overline\partial w)$ in $\tilde{X}_r$. But then,
$\overline\partial \phi(g)=\overline\partial w$ on $X'_r$ and

$$
\int_{X'_r} \overline\partial (\chi\,\psi)\wedge w=\int_{X'_r}
\overline\partial  (\chi\,\psi)\wedge ([w-\phi(g)] +\int_{X'_r}
\overline\partial (\chi\,\psi) \wedge \phi(g)=F+G
$$

\medskip
\noindent Since $w-\phi(g)\in Z^{0,n-1}_{(2)}(X'_r)$ we know that
$F=-(-1)^n \int_{X'_r} \chi\psi\wedge
\overline\partial[w-\phi(g)]=0$. To finish the proof of the
proposition, we need to show that $G:=\int_{X'_r}
\overline\partial(\chi\,\psi)\wedge \phi(g)=-(-1)^n \int_{X'_r}
\chi\psi \wedge \overline\partial \phi(g)$. Using Lemma 4.2, we have
that $\phi(g)\in \text{Dom}(\overline\partial^1)$. Hence there exist
$g_{\nu}\in C^{\infty}_{0,\,(0,n-1)}(\overline{X}_r\setminus\{0\})$ such that
$g_{\nu}\to \phi(g)$. But then, using Stokes' theorem we obtain

\begin{eqnarray*}
\int_{X'_r} \overline\partial (\chi\,\psi)\wedge
\phi(g)&=&\lim_{\nu\to \infty} \int_{X'_r}
\overline\partial(\chi\,\psi)\wedge g_{\nu}=\lim_{\nu\to \infty}
\int_{X'_r} d(\chi\,\psi\,\wedge g_{\nu})-\\
&-& (-1)^n \lim_{\nu \to \infty} \int_{X'_r} \chi\,\psi\wedge
\overline\partial\,g_{\nu}=-(-1)^n \int_{X'_r} \chi\psi \wedge
\overline\partial \phi(g).\\
\end{eqnarray*}

\medskip
\noindent Hence the Proposition is proven. 

\subsection{Towards an understanding of $\text{Im}\,j_*$}  In section 2 of the paper, we defined two maps $j_*,\,\ell_*$ that would be crucial for the proof of Theorem 1.2. The map $j_*:H^{0,n-1}_{(2)}(X'_r)\to H^{n-1}(X'_r,\,\mathcal{O})$ was induced by
the inclusion $j: L^{0,n-1}_{(2)}(X'_r)\to L^{0,n-1}_{2,\,\text{loc}}(X'_r)$ and by Theorem 2.4 we know it   is injective.  In section 2.2 we showed that there
exists a natural map $\ell:
L^{0,n-1}_{2,\,\text{loc}}(\tilde{X}_r,\,\mathcal{O}(D))\to
L^{0,n-1}_{2,\,\text{loc}}(X'_r)$ defined by

$$\ell(g):=(\pi^{-1})^* (g\,).$$

\noindent Clearly $\ell$ commutes with $\overline\partial$ and
induces a map $\ell_*: H^{n-1}(\tilde{X}_r,\,\mathcal{O}(D))\to
H^{n-1}(X'_r,\,\mathcal{O})$ in cohomology.

\medskip
\noindent
We will begin by  obtaining a
characterization for forms $f\in L^{0,n-1}_{2,\text{loc}}(X'_r)\cap
\text{kern}(\overline\partial)$ that arise as $j_*([h])=[f]$ for
some $h\in Z^{0,n-1}_{(2)}(X'_r)$. The lemma below describes some
necessary and  sufficient conditions to address this question.

\medskip
\begin{lemma} Let $f\in L^{0,n-1}_{2,\,\text{loc}}(X'_r)\cap
\text{kern}\,\overline\partial$. Then,

\medskip
\noindent i) If $[f]\in \text{Im}\,j_*$, then $<f,\psi>:=\int_{X'_r}
f\wedge \overline\partial \chi \wedge \psi=0$ when $\psi\in
\text{Kern}(\overline\partial^1)^{n,0}$. Here $\chi$ is a cut-off
function as in section 4.3.

\medskip
\noindent ii) On the other hand, if $\int_{X'_r} |f|^2\,\|z\|^B\,
dV<\infty$ for some $B>0$ large enough and $<f,\,\psi>=0$ when $
\psi\in \text{kern}(\overline\partial^1)^{n,0}$, then $[f]\in
\text{Im}\,j_*$.
\end{lemma}

\medskip
\noindent {\it{Proof.}} i) In section 4.3, we constructed a pairing
$<,\,> : Z^{0,n-1}_{(2)}(X'_r)\times H^{n,0}_{(2)}(X'_r)\to
\mathbb{C}$ described by: for $(f,\,\psi)\in
Z^{0,n-1}_{(2)}(X'_r)\times H^{n,0}_{(2)}(X'_r)$

$$
<f,\,\psi>=\int_{X'_r} f \wedge \overline\partial\chi\,\wedge \psi
$$

\noindent where $\chi\in C^{\infty}_0(X_r)$ such that $\chi=1$ near
$0$. Certainly this pairing can be defined also for forms $f\in
L^{0,\,n-1}_{2,\,\text{loc}}(X'_r)\cap
\text{kern}\,\overline\partial$.

\medskip
\noindent Let $[f]=j_*[h]$ for some $h\in Z^{0,n-1}_{(2)}(X'_r)$.
Then $f=h+\overline\partial u$ for some $u\in
L^{0,\,n-2}_{2,\,\text{loc}}(X'_r)$. In Proposition 4.4, we showed
that whenever $\psi\in \text{kern}\,(\overline\partial^1)^{n,0}$ we
have $<a,\,\psi>=0$ for all $a\in Z^{0,n-1}_{(2)}(X'_r)$. Hence, to
prove i), it suffices to show that $<\overline\partial u,\,\psi>=0$
when $\psi\in \text{kern}\,(\overline\partial^1)^{n,0}$. But this
follows easily from Stokes' theorem as

$$
<\overline\partial u,\psi>=\int_{X'_r} d\,(u\wedge
\overline\partial\chi  \wedge \psi)=0$$

\noindent since the integrand form is compactly supported in $X'_r$.

\medskip
\noindent ii) As $<f,\,\psi>=0$ for all $\psi\in
\text{kern}\,(\overline\partial^1)^{n,0}$, the bounded linear
functional

$$
\lambda: \text{kern}\,(\overline\partial)^{n,0}\to \mathbb{C}
$$

\medskip
\noindent defined by $\lambda(\psi)=<f,\,\psi>=\int_{X'_r}
f\wedge\overline\partial\chi\wedge \psi$ factors to a well-defined
bounded linear functional (still denoted by $\lambda$ for
simplicity)

$$
\lambda:
\dfrac{\text{kern}\,(\overline\partial)^{n,0}}{\text{kern}\,(\overline\partial^1)^{n,0}}\to
\mathbb{C}
$$

\noindent such that $\lambda([\psi])=<f,\psi>$.

\medskip
\noindent In Proposition 4.5, we saw that

$$
\left(\dfrac{\text{kern}\,(\overline\partial)^{n,0}}{\text{kern}\,(\overline\partial^1)^{n,0}}\right)'
\cong \dfrac{H^{0,n-1}_{(2)}(X'_r)}{\phi_*(H^{0,n-1}(\tilde{X}_r))}.
$$

\noindent Hence there exists a $g\in Z^{0,n-1}_{(2)}(X'_r)$ such
that $\lambda([\psi])=<g,\,\psi>$ for all $\psi\in
H^{n,0}_{(2)}(X'_r)$., i.e.

\begin{eqnarray}\label{eq:mc}
\int_{X'_r} f\wedge \overline\partial \chi \wedge \psi=\int_{X'_r}
g\wedge \overline\partial \chi\,\wedge \psi
\end{eqnarray}

\noindent for all $\psi\in H^{n,0}_{(2)}(X'_r)$.

\medskip
\noindent Arguing now verbatim as in section 4.2, condition
$(\ref{eq:mc})$ will guarantee the existence of a $w\in
L^{0,n-1}_{(2)}(\tilde{X}_r)$, compactly supported in $\tilde{X}_r,$
such that $\overline\partial w=\pi^*\,\left(\overline\partial
\chi\wedge(f-g)\right)$ or equivalently the existence of a
$u=\phi(w)\in \text{Im}\,\phi$, compactly supported in $X_r,$
satisfying $\overline\partial u=\overline\partial \chi \wedge (f-g)$
on $X'_r$. Then, as in section 4.2, we can split $f-g=f_1+f_2$,
where

\begin{eqnarray*}
f_1&=&\chi(f-g)-u,\,\;\;\text{compactly\;\;supported\;\;in\;\;}\;X_r,\\
f_2&=& (1-\chi)\,(f-g)+u\in L^{0,n-1}_{(2)}(X'_r)\\
\end{eqnarray*}

\medskip
\noindent and $f_1,\,f_2$ are $\overline\partial$-closed. Moreover,
we have $\int_{X'_r} |f_1|^2\;\|z\|^B\,dV<\infty$.

\medskip
\noindent We shall now recall a result about weighted
$L^2$-estimates for solutions to $\overline\partial$-closed forms defined on
$\text{Reg}\,\Omega$, compactly supported in $\Omega$, where
$\Omega$ Stein relatively compact subdomain of a Stein space and
$A=\text{Sing}\,X$:

\begin{theorem}(Theorem 5.3 in \cite{OV2}{\footnotemark\footnotetext{Theorem 5.3 in
\cite{OV2} works more generally for $A$ complex analytic sets
containing the singular locus of $X$.}}) Let $f$ be a $(p,q)$ form
defined on $\text{Reg}\,\Omega$ and $\overline\partial$-closed there
with $0<q<n$, compactly supported in $\Omega$ and such that
$\int_{\text{Reg}\,\Omega} |f|^2\, d_A^{N_0} dV<\infty$ for some
$N_0\ge 0$. Then there exists a solution $u$ to $\overline\partial
u=f$ on $\text{Reg}\,\Omega$ satisfying $\text{supp}_X\, u\Subset
\Omega$ and such that

$$\int_{\text{Reg}\,\Omega} |u|^2\, d_A^{N} dV\le
C\,\int_{\text{Reg}\,\Omega} |f|^2\, d_A^{N_0} dV
$$

\noindent where $N$ is a positive integer that depends on $N_0$ and
$\Omega$ and $C$ is a positive constant that depends on
$N_0,\,N,\,\Omega$\; and\;\;$\text{supp}\,f$. Here $d_A$ denotes the
distance function to $A$.

\end{theorem}

\medskip
\noindent Using the above theorem, we know there exists a $v\in
L^{0,n-2}_{2,\text{loc}}(X'_r)$ such that $\overline\partial v=f_1$.
Therefore we have

$$f=\left((1-\chi)\,f+\chi\,g+u\right)+\overline\partial v.$$

\noindent Hence $[f]=j_*[\,(1-\chi)\,f+\chi\,g+u\,]$. Q.E.D. 

\medskip
\noindent \bigskip \noindent In \cite{Rup3}, Ruppenthal identified more or less
the $\text{kern}(\overline\partial^1)^{n,0}$ in terms of resolution
data. Using the notation of Theorem 1.2 in the introduction we have:

\begin{lemma} (Lemma 6.2 in \cite{Rup3})
$\text{kern}(\overline\partial^1)^{n,0}=(\pi^{-1})^*
\left(L^{n,0}_{(2)}(\tilde{X}_r)\cap
\Gamma(\tilde{X}_r,\,\mathcal{K}_{\tilde{X}_r}\otimes
\mathcal{O}(|Z|-Z))\right)$.
\end{lemma}

\medskip
\noindent {\bf{Remark:}} Lemma 6.2 in \cite{Rup3} only states that
$\Gamma(\tilde{X}_r,\,\mathcal{K}_{\tilde{X}_r}\otimes
\mathcal{O}(|Z|-Z))=\text{kern}(\overline\partial_{s,\,loc})\cap
\Gamma(\tilde{X}_r,\,\mathcal{F}^{n,0}_{\gamma,\,E})$, where
$\overline\partial_{s,loc}$ is defined as follows: Let $f\in
L^{p,q}_{\gamma,\,loc}(\tilde{X}_r)$, where $\gamma$ is the
``pseudometric'' from section 2. We say that $f\in
\text{Dom}(\overline\partial_{s,\,loc})(\tilde{X}_r)$ if
$\overline\partial f\in L^{p,q+1}_{\gamma,\,loc}(\tilde{X}_r)$ and
there exist a sequence of smooth forms $f_j$ compactly supported
away from $E$ such that $f_j\to f$ in the graph norm in
$L^{p,*}_{\gamma}(K)$ for any $K$ compact subset of $\tilde{X}_r$.
We write in this case
$\overline\partial_{s,\,loc}f=\overline\partial f$ (where the right
hand side is taken in the weak sense). Also
$\mathcal{F}^{n,0}_{\gamma,\,E}(\tilde{X}_r):=L^{n,0}_{\gamma,\,loc}(\tilde{X}_r)\cap
\text{Dom}(\overline\partial_{s,\,loc})(\tilde{X}_r)$.

\medskip
\noindent{\it{Proof of Lemma 4.8.}} It follows immediately from
Lemma 6.2 in \cite{Rup3}  that
$\text{kern}(\overline\partial^1)^{n,0}\subset (\pi^{-1})^*
(L^{n,0}_{(2)}(\tilde{X}_r)\cap $

\noindent $\Gamma(\tilde{X}_r, \,\mathcal{K}_{\tilde{X}_r}\otimes
\mathcal{O}(|Z|-Z)))$. To prove the reverse inclusion, assuming that
$L^{n,0}_{(2)}(\tilde{X}_r)\cap
\Gamma(\tilde{X}_r,\,\mathcal{K}_{\tilde{X}_r}\otimes
\mathcal{O}(|Z|-Z))$ is non-trivial, we proceed as follows: By
Ruppenthal's result we know that given $f\in
\Gamma(\tilde{X}_r,\,\mathcal{K}_{\tilde{X}_r}\otimes
\mathcal{O}(|Z|-Z))$, there exists a sequence of $f_j$ smooth
compactly supported away from $E$ such that $f_j\to f$ in the graph
norm in $L^{n,*}_{\gamma}(K)$ for every compact subset $K$ of
$\tilde{X}_r$. Choose a cut-off function $c\in
C^{\infty}_{0}(\tilde{X}_r)$ such that $c=1$ near $E$. Then
$c\,f_{j}\to c\,f$ in graph norm in $L^{n,*}_{\gamma}(\tilde{X}_r)$
and the same holds true for their push-forward. The push-forward of
$(1-c)\,f$ is easily approximated in graph norm by smooth forms
supported away from $0$ (as $f$ is assumed to be now in
$L^{n,0}_{\gamma}(\tilde{X}_r)=L^{n,0}_{(2)}(\tilde{X}_r))$. Hence
$(\pi^{-1})^*\,f\in \text{kern}(\overline\partial^1)^{n,0}$.

\medskip
\noindent
\subsubsection{\bf{An alternative description of $\text{Im}\,j_*$}}  The second key ingredient in the proof of Theorem 1.2 is the realization that $\text{Im}\,j_*=\text{Im}\,l_*$ or equivalently

\begin{lemma} The map $\ell_*:
H^{n-1}(\tilde{X}_r,\,\mathcal{O}(D))\to
H^{n-1}(X'_r,\,\mathcal{O})$ is surjective on $\text{Im}\,j_*$.
\end{lemma}

\medskip\noindent {\it{Proof.}}  We need to show that i) $\text{Im}\,\ell_*\subset 
\text{Im}\,j_*$ and ii) $\text{Im}\,j_* \subset \text{Im}\ell_*$. To prove i)  it suffices to
show  for any $g\in
L^{0,n-1}_{2,\,\text{loc}}(\tilde{X}_r,\,\mathcal{O}(D))\cap
\text{kern}\,(\overline\partial),$ that $\ell(g)$ satisfies the
conditions of Lemma 4.6 ii), i.e.

\begin{eqnarray*}
&\alpha)&\; \int_{X'_r} |\ell(g)|^2 \|z\|^B\, dV<\infty,\\
&\beta)&\; <\ell(g),\,\psi>=0,\;\;\text{for\;all}\;\psi\in
\text{kern}(\overline\partial^1)^{n,0}
\\
\end{eqnarray*}

\medskip
\noindent Property $\alpha)$ follows easily, for some $B>0$
sufficiently large, by the estimates in section 3 of an earlier
paper of ours, see Lemma 3.1 in \cite{FOV1}. There, we compared
weighted $L^2$-norms between forms and their pull-backs under
resolution of singularities maps.

\medskip
\noindent It remains to prove $\beta)$. When $\psi\in
\text{kern}(\overline\partial^1)^{n,0}$  let
$\tilde{\psi}:=\pi^*\,\psi$. Lemma 4.8 yields immediately that   $\tilde{\psi}\in L^{n,0}_{(2)}(\tilde{X}_r)\cap
\Gamma(\tilde{X}_r,\,\mathcal{K}\otimes \mathcal{O}(-D))$. Let also
$\tilde{\chi}:=\chi\circ \pi$. Then

$$
<\ell(g),\,\psi>=-\int_{\tilde{X}_r} \overline\partial
\tilde{\chi}\wedge g\wedge \tilde{\psi}.
$$

\medskip
\noindent But then $g\wedge \tilde{\psi}\in
L^{n,\,n-1}_{2,\text{loc}}(\tilde{X}_r)$ is
$\overline\partial$-closed outside $E$; thus it extends as a
$\overline\partial$-closed form $b$ in $\tilde{X}_r$. Hence

$$<\ell(g),\,\psi>=-\int_{\tilde{X}_r} \overline\partial
\tilde{\chi}\wedge b=-\int_{\tilde{X}_r} d(\tilde{\chi}\,b)=0
$$

\noindent by Stokes' theorem, since $\tilde{\chi}\,b$ has compact
support in $\tilde{X}_r$. 

\medskip
\noindent To prove ii)\; i.e. that $\text{Im}\,j_*\subset \text{Im}\,\ell_*$ 
we will use a ``twisted'' version of arguments that appeared in sections
4.2-4.4. Let $f\in Z^{0,n-1}_{(2)}(X'_r)$ and let $\chi\in C^{\infty}_0(X_r)$
such that $\chi=1$ near $0$ and $\text{supp}\,\chi\subset X_{\rho}$ for some 
$0<\rho<r$. Let $\tilde{f}:=\pi^*\,f$ and $\tilde{\chi}:=\chi\circ \pi$. Then for all $\psi\in \text{kern}(\overline\partial^1)^{n,0}_{X'_{\rho}}$
we have that $<f,\,\psi>_{X'_{\rho}}=\int_{X'_{\rho}} f\wedge \overline\partial \chi \wedge \psi=0$. Using 
Lemma 4.8 this implies that 

\begin{eqnarray}\label{eq:altmom}
 \int_{\tilde{X_{\rho}}} \tilde{f}\wedge \overline\partial \tilde{\chi}\wedge \tilde{\psi}=0
\end{eqnarray}

\noindent for all $\tilde{\psi}\in \mathcal{L}^{n,0}(\tilde{X}_{\rho},\,\mathcal{O}(-D))\cap \text{kern}\,\overline\partial$ and where $D:=Z-|Z|$ is as in Theorem 1.2. 

\medskip
\noindent We may consider $A:=\tilde{f}\wedge \overline\partial \tilde{\chi}\otimes s$ as an element in $L^{0,n}_{(2)}(\tilde{X}_{\rho},\,L_D)$ with $\text{supp}\,u\subset \tilde{X}_{\rho}$ and let $B:=\tilde{\psi}\otimes s^{-1}$\\
\noindent $\in L^{n,0}_{(2)}(U,\,L_{-D})$. Then we can rewrite $(\ref{eq:altmom})$   as 

\begin{eqnarray}\label{eq:gmom}
\int_U A\wedge B=0
\end{eqnarray}

\medskip
\noindent  Now, the generalized moment condition $(\ref{eq:gmom})$ will permit us to solve the equation $\overline\partial_D F =A=\tilde{f}\wedge \overline\partial \tilde{\chi}\otimes s$ with $F\in L^{0,n-1}_{(2)}(\tilde{X}_r,\,L_D)$ and $\text{supp}\,F\subset \overline{\tilde{X}_{\rho}}\subset \tilde{X}_r$. This is a consequence  of the following $L^2$-Cauchy problem: 

\begin{proposition} Let $U\subset \subset \tilde{X}_r$ be an open neighborhood of $E$ with smooth strongly
 pseudoconvex boundary, and let $h\in L^{0,n}_{(2)}(\tilde{X}_r,\,L_D),\;\text{supp}\,h\subset \subset U$. If

\begin{eqnarray}\label{eq:mom}
\int_{U} h\wedge a=0\;\;\;\text{for\;
all}\;\;a \in L^{n,0}_{(2)}(U, L_{-D})
\end{eqnarray}

\noindent then there exists a solution $v\in L^{0,n-1}_{(2)}(\tilde{X}_r,\,L_D)$ satisfying $\overline\partial v=h$ with $\text{supp}\,v\subset \overline{U}\subset\subset  \tilde{X}_r$.
\end{proposition} 

\smallskip\noindent {\it{Proof.}} Since $U$ is a smoothly bounded strongly pseudoconvex domain in a complex manifold
we know from Lemma 2.2 of section 2, that the $\text{Range}(\overline\partial_{-D})$ is closed in $L^{n,1}_{(2)} (U,\,L_{-D})$. Hence we have the following strong decomposition

$$L^{n,0}_{(2)}(U,\,L_{-D})=\text{Rang}(\overline\partial^*_{-D}) \oplus H^{n,0}_{(2)}(U,L_{-D})$$

\noindent where $H^{n,0}_{(2)}(U,L_{-D}):=\text{kern}(\overline\partial_{-D})\cap L^{n,0}_{(2)}(U,L_{-D})$.

\medskip\noindent Now, $(\ref{eq:mom})$ implies that $\overline{*}_{D} h$ is orthogonal to $H^{n,0}_{(2)}(U,\,L_{-D}))$, hence it 
belongs to the range of $\overline\partial^*_{-D}$, i.e.  there exists an element $a\in \text{Dom}(\overline\partial^*_{-D})$ such that 
$\overline\partial^*_{-D}a=\overline{*}h$.   By Proposition 1 in \cite{CS}, we know that $\overline{*}_{-D}a\in \text{Dom}(\overline\partial_{D,min})$ and $\overline\partial^*_{-D}=-\overline{*}_{D}\,(\overline\partial_{D,min})\,\overline{*}_{-D}$. Hence 
there exist a sequence of compactly supported sections  $a_n\in D^{0,n-1}(U,L_{D})$ such that $a_n\to \overline{*}_{-D} a$ and $\overline\partial_D a_n\to -(-1)^n \overline{*}_{-D}\,\overline{*}_D\,h=-h$ in $L^{\bullet,\bullet}_{(2)}(U,\,L_D)$. 

\medskip
\noindent Set $\hat{a}:=(-\overline{*}_{-D}a)^0$, i.e. the trivial extension by zero outside $U$. Then we claim that 
$\overline\partial\hat{a}=h$ in $\tilde{X}_r$. Indeed, take $\psi\in D^{0,n}(\tilde{X}_r,\,L_D)$ and let us look at 

\begin{eqnarray*}
(\hat{a},\vartheta_D\,\psi)_{\tilde{X}_r}&=&(-\overline{*}_{-D}\,a,\,\vartheta_D\,\psi)_U=-\text{lim}(a_n, \vartheta_D \psi)_U\\
&-&\text{lim}(\overline\partial_D a_n,\psi)_U=(h,\,\psi)_U=(h,\psi)_{\tilde{X}_r}.\\
\end{eqnarray*}

\noindent Here we used the fact that $a_n$ are compactly supported in $U$ in order to perform integration by parts in the second line, and that $h$ is compactly supported in $U$.

%Hence from the above decomposition we can write (taken as $E:=L_{-D}$ and simplifying $\overline{*}_{L_D}:=\overline{*}_D$)

%$$
%\overline{*}_D \,h=(\overline\partial\, \overline\partial^*_{-D}+\overline\partial^*_{-D}\,\overline\partial)\, N\,\overline{*}_D\,h.
%$$

%\medskip
%\noindent Set $v=\overline{*}_{-D}\overline\partial N\,\overline{*}_{D}\,h$ in $U$ and let $v=0$  in $
%\tilde{X}_r\setminus U$. Then one can easily show that $\overline\partial v= h$ in $\tilde{X}_r$ and $\text{supp}\,v\subset 
%\overline{U}\subset\subset \tilde{X}_r$.  Q.E.D.

\medskip
\noindent Using the above proposition  for $U=\tilde{X}_{\rho}$,  we obtain  a solution $w$ to 
$\overline\partial w=u=\tilde{f}\wedge \overline\partial \tilde{\chi}$ with $w\in L^{0,n-1}_{(2)}(\tilde{X}_r,\,\mathcal{O}(D))$ and $\text{supp}\,w\subset \tilde{X}_{r}$. 

\medskip
\noindent Then, we can write $\tilde{f}=(\tilde{\chi}\,\tilde{f}-w)+(w+(1-\tilde{\chi})\,\tilde{f})=:g_1+g$. 
Let $h:=(\pi^{-1})^{*}\, g_1$. Then $\text{supp}\,h\subset\subset X_{r},\;\; \overline\partial h=0$ on $X'_r$
and $\int_{X'_r} \|z\|^B |h|^2\, dV<\infty$ for some $B>0$ sufficiently large. Then, by Theorem 4.7 we know that 
there exists a solution $t\in L^{0,n-2}_{2,\,\text{loc}}(X'_r)$ such that $\overline\partial t=h$. Hence we can write 

$$f=\ell(g)+\overline\partial\,t.$$

\noindent Hence we have $j_*([f])=\ell_*([g])$ and thus ii) is proven. 

\medskip
\noindent  Let $[f_1],\cdots [f_m]$ be a basis of $H^{0,n-1}_{(2)}(X'_r)$.  Then we can 
define a map $S: H^{0,n-1}_{(2)}(X'_r)\to H^{n-1}(\tilde{X}_r,\,\mathcal{O}(D))$ such that $S(\sum c_j [f_j])=\sum c_j[g_j]$, where 
$[g_j] $ satisfy $j_*[f_j]=\ell_*[g_j]$.

\bigskip \noindent

\subsection{Proof of Theorem 1.2} The sheaf inclusion 
$m: \mathcal{O}_{\tilde{X}_r}\to \mathcal{O}_{\tilde{X}_r}(D)$
induces a commutative diagram between long exact local cohomology
sequences

\xymatrix{& H^{n-2}(\tilde{X}_r\setminus
E,\,\mathcal{O})\ar[d]_{m_*}^{\cong} \ar[r]^{\delta}
&H^{n-1}_E(\tilde{X}_r,\,\mathcal{O})\ar[d]_{m_{*,\,E}}\ar[r]^{k_*}
&H^{n-1}(\tilde{X}_r,\,\mathcal{O}) \ar[d]_{m_*} \ar[r]^{r_*}
&H^{n-1}(\tilde{X}_r\setminus E,\,\mathcal{O})\ar[d]_{m_*}^{\cong}&^{}\\
& H^{n-2}(\tilde{X}_r\setminus E,\,\mathcal{O}(D)) \ar[r]^{\delta} &
H^{n-1}_E(\tilde{X}_r,\,\mathcal{O}(D))\ar[r]^{k_*}
&H^{n-1}(\tilde{X}_r,\,\mathcal{O}(D))\ar[r]^{r_*}
&H^{n-1}(\tilde{X}_r\setminus E,\,\mathcal{O}(D)).\\
}

\medskip
\noindent By Karras' result we know that
$H^{n-1}_E(\tilde{X}_r,\,\mathcal{O})=0$. Taking into account this
and the commutativity of the above diagram (in particular of the
left square) we obtain the following exact sequence:

\begin{eqnarray}\label{eq:lastle}
0\to H^{n-1}_E(\tilde{X}_r,\,\mathcal{O}(D))\overset{k_*}
\longrightarrow
H^{n-1}(\tilde{X}_r,\,\mathcal{O}(D))\overset{r_*}\longrightarrow
H^{n-1}(\tilde{X}_r\setminus E,\,\mathcal{O}(D)).
\end{eqnarray}

\medskip
\noindent
\medskip
\noindent Using this information
we construct the following diagram:

\xymatrix{ 0\ar[r]
&H^{n-1}_E(\tilde{X}_r,\,\mathcal{O}(D))\ar[r]^{k_*}
&H^{n-1}(\tilde{X}_r,\,\mathcal{O}(D))\ar@{-->>}[dl]_{T}
\ar@{-->>}[d]^{\ell_*} \ar[dr]_{\ell_*} \ar[r]^{r_*}
&H^{n-1}(\tilde{X}_r\setminus
E,\,\mathcal{O}(D))&^{}\\
^{} &  H^{0,n-1}_{(2)}(X'_r) \ar[r]^{j_*}_{\cong}
&(\text{Im}\,j_*)\ar@{^(->}[r] &H^{n-1}(X'_r,\,\mathcal{O})\ar[r]^{\pi^*}_{\cong} &H^{n-1}(\tilde{X}_r\setminus E,\,\mathcal{O})\ar[ul]^{m_*}_{\cong}.\\
}

\medskip
\noindent The top row is the exact sequence  from
$(\ref{eq:lastle})$. From the right rectangle of the diagram, we
observe that $r_*=m_*\circ (\pi^*)\circ \ell_*$.

\medskip
\noindent We shall show in a moment that

\begin{lemma}  i) The natural map $k_*:
H^{n-1}_E(\tilde{X}_r,\,\mathcal{O}(D))\to
H^{n-1}(\tilde{X}_r,\,\mathcal{O}(D))$ is injective with
$\text{Im}\,k_*=\text{kern}\,\ell_*$.

\medskip
\noindent ii) The map $j_*^{-1}\circ \ell_* :
H^{n-1}(\tilde{X}_r,\,\mathcal{O}(D))\to
H^{0,n-1}_{(2)}(\tilde{X}_r)$ is a surjective map $T$ as in Theorem

1.2.
\end{lemma}

\medskip
\noindent{\it{Proof.}} i) The injectivity of $k_*$ follows from the
exactness of $(\ref{eq:lastle})$. Also, from the exactness of
$(\ref{eq:lastle})$ we have that $\text{Im}\,k_*=\text{kern}\,r_*$.
Due to the commutativity of the right rectangle of the above
diagram, we see that $\text{kern}\,r_*=\text{kern}\,\ell_*$.

\medskip
\noindent ii) As $j_*$ is an isomorphism between
$H^{0,n-1}_{(2)}(X'_r)$ and $\text{Im}\,j_*$, we can define the
map $T:=j_*^{-1}\circ \ell_*:
H^{n-1}(\tilde{X}_r,\,\mathcal{O}(D))\to H^{0,n-1}_{(2)}(X'_r)$.
Clearly $\text{kern}\,T=\text{kern}\,\ell_*=\text{Im}\,k_*\cong
H^{n-1}_E(\tilde{X}_r,\,\mathcal{O}(D))$. The surjectivity of $T$ 
follows from the fact that $T\circ S=Id$ on $H^{0,n-1}_{(2)}(X'_r)$. 

\medskip
\noindent{\bf{Remark 4.5.1}} When $q< n-1$ the map $j_*: H^{0,q}_{(2)}(X'_r)\to H^q(X'_r,\,\mathcal{O})$ is an isomorphism. Arguing in a similar manner as in 
section 4.5 we obtain  the following 
short exact sequence of sheaves for each $q\le n-2$: 

$$0\to H^q_E(\tilde{X}_r,\,\mathcal{O}(D))\to H^q(\tilde{X}_r,\,\mathcal{O}(D))\to H^{0,q}_{(2)}(X'_r)\to H^{q+1}_E(\tilde{X}_r,\,\mathcal{O}(D))$$

\noindent where the $H^{0,q}_{(2)}(X'_r)$ entry  appears due to the fact that $H^{q}(\tilde{X}_r\setminus E,\,\mathcal{O})\cong H^q(X'_r,\,\mathcal{O})\overset{{j^{-1}_*}}
\cong H^{0,q}_{(2)}(X'_r)$.

\medskip
\noindent As a consequence of the above sequence and in the special case where $-D$ is locally semi-positive with respect to $X$ we obtain 
(via Theorem 3.1 and Takegoshi's vanishing theorem) for all $q$ with 
$0\le q\le n-2$ that $H^{0,q}_{(2)}(X'_r)\cong H^q(\tilde{X}_r, \,\mathcal{O}(D))$. The isomorphism when $q=n-1$ in this case has already been
observed in the introduction  as a consequence of Theorem 1.2. Hence we can recover Ruppenthal's Theorem  7.1 from \cite{Rup3} for all 
$q\le n-1$. 

\medskip
\noindent
\section{Proofs of Theorems 1.3 and 1.4} 

\medskip
\noindent \subsection{Proof of Theorem 1.3.}  We choose neighborhoods $\{V_j\}_{j=1}^m$ of $\{a_j\}_{j=1}^m$ with 
$V_j\subset\subset X$ and such that for all $j=1,\cdots m, \;V_j \cong \hat{X} ^j_r\subset B(0,R)\subset \mathbb{C}^{N_j}$,  where $\hat{X} ^j_r$ are subvarieties with $0$ as an isolated singular point.  Assume $\overline{V_i}\cap \overline{V}_j=\emptyset$, if $i\neq j$. Set 
$V:=\cup_{j=1}^m V_j,\;\tilde{V}=\pi^{-1} (V)$. Choose a partition of unity $\chi_0,\,\chi_1,\cdots, \chi_m$ with $\text{supp}\,\chi_0\subset \overline{X}\setminus \Sigma$ and $\text{supp}\chi_j\subset V_j$ if $j>0$, Thus $\chi_j=1$ near $a_j$ and $\chi_0=1$ 
near $\overline{X}\setminus V$. Let $\phi_*: H^{0,q}_{(2)}(\tilde{X},\,\mathcal{O})\to H^{0,q}_{(2)}(X')$ be the map sending 
$[f]$ to $[(\pi^{-1})^*f]$. We need to show that $\phi_*$ is bijective for $1\le q\le n-2$.

\medskip
\noindent  We show first surjectivity. Let $[f]\in H^{0,q}_{(2)}(X')$. By theorem 1.1, we know that $f_{\upharpoonright_{V_i}}=\phi(g_i)+\overline\partial u_i$ for $i=1,\cdots,m$
where $g_i\in Z^{0,q}_{(2)}(\tilde{\hat{X}}^i_r)$ and $u_i\in L^{0,q-1}_{(2)}(\hat{X}^i_r)$. Set  $g:=\pi^*f$ on $\tilde{X}\setminus \tilde{V}$
and $g:=g_i+\overline\partial \pi^*(\chi_o\,u_i)$, on $V_i$. Then $g\in Z^{0,q}_{(2)}(\tilde{X})$ and $f-\phi(g)=\sum_1^m \overline\partial (\chi_i\,u_i)^{o}$ where by ${k}^o$ we mean trivial extension of a function $k$ by zero outside $V_i$. Then $[f]=\phi_*[g]$.

\medskip
\noindent To show injectivity, we let $g\in Z^{0,q}_{(2)}(\tilde{X})$ and assume that $\phi(g)=\overline\partial u$ for some 
$u\in L^{0,q-1}_{(2)}(X')$. Write $g:=\overline\partial (\pi^*(\chi_0\,u))+\sum_1^m g_i$, where $g_i=\tilde{\chi}_i\,g+\overline\partial\tilde{\chi}_i\wedge \pi^* u$. We have $g_i\in L^{0,q}_{(2)}(\tilde{V_i})$ with $\overline\partial g_i=0$ and $\text{supp}\,g_i$ compact in $\tilde{V}_i$.  Hence $[g_i]\in H^{q}_c(\tilde{V}_i,\,\mathcal{O})$. If $A_i=\pi^{-1}(a_i)$ is the exceptional set of the desingularization $\pi: \tilde{V}_i \to V_i$, then  by Karras' result we have that $H^q_c(\tilde{V}_i,\,\mathcal{O})=0$. Hence there exists $v_i\in L^{0,q-1}_{(2)}(\tilde{V}_i)$, compactly supported in $\tilde{V}_i$ such that $\overline\partial v_i=g_i$. Then we set  $v:=(\pi)^*(\chi_0\,u)+\sum v_i^{o}$. We can easily check that 
$v\in L^{0,\,q-1}_{(2)}(\tilde{X})$ and $\overline\partial v=g$; hence $[g]=0$.  

\medskip
\noindent \subsection{Proof of Theorem 1.4.}  In the proof of Theorem 1.2,  the map $j_*: H^{p,q}_{(2)}(X'_r)\to 
H^{q}_{(2)}(X'_r, \Omega^p)$ which was induced by the inclusion  $j:L^{p,q}_{(2)}(X'_r)\to L^{p,q}_{2,\,\text{loc}}(X'_r)$
played a crucial role. For the situation we consider in Theorem 1.4, we need to introduce some auxiliary spaces and a modified map
$j'_*$. More precisely, let us set 

$$
^{\prime} L^{p,q}_{(2)}(X'):=\{f\in L^{p,q}_{2,\,\text{loc}}(X')\;|\;\; f\in L^{p,q}_{(2)}(X\setminus V)\}
$$

\medskip
\noindent Let $^{\prime} \mathcal{L}^{p,q}(X'):=\{f\in\;\;    ^{\prime} L^{p,q}_{(2)}(X'),\;\;\overline\partial f\in \;\;  ^{\prime} L^{p,q+1}_{(2)}(X')\}$ and let
$^{\prime} H^{p,q}(X')$ denote the cohomology of the complex $\left(^{\prime} \mathcal{L}^{p,\bullet,},\,\overline\partial\right)$, where $\overline\partial$ is taken with respect to the open subsets in $X'$. Let  us consider the inclusion map $j': L^{p,q}_{(2)}(X') \to \;
^{\prime} L^{p,q}_{(2)}(X')$.  Then we have: 

\begin{proposition} For $p+q\le n-1,\,q>0$ the map $j'_*: H^{p,q}_{(2)}(X')\longrightarrow\;^{\prime}H^{p,q}(X')$ is injective.
\end{proposition}

\medskip
\noindent{\it{Proof.}} Let $f\in Z^{p,q}_{(2)}(X')$ and assume that $j'_*([f])=0$, i.e. $\overline\partial u=f$ for some 
$u\in\;^{\prime}L^{p,q-1}_{(2)}(X')$. Using the partition of unity $\{\chi_i\}_{i=0}^m$, we can rewrite $f$ as 

$$f=\overline\partial u=\overline\partial\, (\sum_{i=0}^m \chi_i\,u)=\overline\partial (\chi_0\,u)+\sum_{i=1}^m (\chi_i\,f+\overline\partial
\chi_i\,\wedge u).$$ 

\medskip 
\noindent Now the forms  $g_i:=\chi_i\,f+\overline\partial\chi_i\,\wedge u\in L^{p,q}_{(2)}(V_i\setminus\{a_i\})$ are $\overline\partial$-
closed there and $\text{supp}\,g_i\subset\subset V_i$; hence by Proposition 3.1  in \cite{FOV}, we know that there exists a  
$v_i\in L^{p,q-1}_{(2)}(V_i\setminus\{a_i\})$, with compact support in $V_i$ such that $\overline\partial v_i=g_i$ for $i=1,\cdots,m$. 
Setting $v:=\chi_0\,u+\sum_{i=1}^m\,v_i\in L^{p,q-1}_{(2)}(X')$ we have  $\overline\partial v=f$. \;\;\;\;Q.E.D.

\medskip
\noindent{\bf{Remark 5.2.1:}} Using a similar argument one can further show that the map $j'_*: H^{p,q}_{(2)}(X')\longrightarrow\; ^{\prime}H^{p,q}(X')$ is bijective for $p+q\le n-2,\;q>0$. 

\medskip 
\noindent Let us consider the map $\ell^{\prime}: L^{0,n-1}_{(2)}(\tilde{X},\,\mathcal{O}(D))\to\;^{\prime}L^{0,n-1}_{(2)}(X')$ which sends 
$g\to (\pi^{-1})^*g$ and let $\ell^{\prime}_*: H^{0,n-1}_{(2)}(\tilde{X},\,\mathcal{O}(D))\to\;^{\prime}H^{0,n-1}(X')$ be the corresponding map in cohomology. The first step in the proof of Theorem 1.4 is to show that $\text{Im}\,(j'_*)=\text{Im}\,(\ell^{\prime}_*)$. 

\medskip
\noindent We shall show first that $\text{Im}\,(\ell^{\prime}_*)\subset \text{Im}\,(j'_*)$. Let us consider an element 
$g\in Z^{0,n-1}_{(2)}(\tilde{X},\,\mathcal{O}(D)):=L^{0,n-1}_{(2)}(\tilde{X},\,\\$
\noindent
$\mathcal{O}(D))\cap \text{kern}(\overline\partial)$. 
By Lemma 4.9, we know that on $V_i$  we have $\ell'(g)=f_i+\overline\partial u_i$ where $f_i\in Z^{0,n-1}_{(2)}(V_i\setminus a_i)$ and $u_i\in L^{0,n-2}_{2,\text{loc}}(V_i\setminus a_i)$. Set $f:=\ell'(g)$ on $X\setminus V$ and $f:=f_i+\overline\partial(\chi_0\,u_i)$ on $V_i\setminus a_i$. Then $f\in L^{0,n-1}_{(2)}(X')$, is well-defined and $\overline\partial$-closed and $\ell'(g)-f=\sum_{i=1}^m\,\overline\partial\,(\chi_i\,u_i)^o$. Hence $j'_*([f])=\ell^{\prime}_*[g]$. 

\medskip
\noindent To show the other direction, we consider an element $f\in Z^{0,n-1}_{(2)}(X')$. By Lemma 4.9, we have on each $V_i;\;
f=\ell^{\prime}(g_i)+\overline\partial u_i$, where $g_i\in Z^{0,n-1}_{(2)}(\tilde{V}_i,\,\mathcal{O}(D))$ and $u_i\in L^{0,n-2}_{2,\text{loc}}(V_i\setminus a_i)$. Set $g:=\pi^*f$ on $\tilde{X}\setminus \tilde{V},$ and $g:=g_i+\overline\partial \pi^*(\chi_0\,u_i)$ 
on $\tilde{V}_i$ for $i=1,\cdots,m$. Then $g\in L^{0,n-1}_{(2)}(\tilde{X},\,\mathcal{O}(D)),$ is well-defined with $\overline\partial g=0$. 
Then $f=\ell^{\prime}(g)+\sum_{i=1}^m \overline\partial (\chi_i\,u_i)^o$; hence $j'_*[f]=\ell'_*[g]$. 

\medskip
\noindent Then we can consider the operator $\tilde{T}: H^{0,n-1}_{(2)}(\tilde{X},\,\mathcal{O}(D))\to H^{0,n-1}_{(2)}(X')$ defined by
$\tilde{T}:={j'_*}^{-1}\circ \ell'_*$; clearly $\tilde{T}$ is surjective. It remains to show that the kernel of $\tilde{T}$ is naturally isomorphic to $H^{n-1}_E(\tilde{X},\mathcal{O}(D))$.  As in the local case, we have the following short exact sequence

$$
0\to H^{n-1}_E(\tilde{X},\,\mathcal{O}(D))\to H^{n-1}(\tilde{X},\,\mathcal{O}(D))\to H^{n-1}(\tilde{X}\setminus E,\,\mathcal{O}(D))
$$

\medskip
\noindent By Karras' result we know that $H^{n-1}_E(\tilde{X},\,\mathcal{O}(D))\cong H^{n-1}_c(\tilde{V},\,\mathcal{O}(D))$ and we have the following commutative diagram: 

\xymatrix{ 0\ar[r]
&H^{n-1}_E(\tilde{X},\,\mathcal{O}(D))\ar[d]^{\cong} \ar[r]
&H^{n-1}(\tilde{X},\,\mathcal{O}(D)) \ar[r]^{r_*}
&H^{n-1}(\tilde{X}_r\setminus
E,\,\mathcal{O}(D))&^{}\\
^{} &  H^{n-1}_{c}(\tilde{V},\,\mathcal{O}(D))\ar@{-->}[ur]  \ar[r]^{i'_*} & H^{0,n-1}_{(2)}(\tilde{X},\,\mathcal{O}(D))\ar[u]^{i''_*} \ar[r]^{\ell'_*} &^{\prime}H^{0,n-1}(X').\\
}

\medskip
\noindent 
From the commutativity of the left upper triangular part of the diagram we can conclude that the oblique map from $H^{n-1}_c(\tilde{V},\mathcal{O}(D))\to H^{n-1}(\tilde{X},\,\mathcal{O}(D))$ is injective. Similarly, from the commutativity of the bottom left triangular 
part of the diagram we can conclude that $i'_*$ is injective. Now, by the definition of $\tilde{T}$ we know that $\text{kern}\,\tilde{T}=\text{kern}\,\ell'_*$. If we could show that the 
$\text{kern}\,\ell'_*=\text{Im}\,i'_*$, then   we could finish the proof of Theorem 1.4, since the map $i'_*$ is injective and hence $\text{kern}\,\tilde{T}\cong H^{n-1}_E(\tilde{X},\,\mathcal{O}(D))$. 

\medskip
\noindent We first observe that $\ell'_*\circ i'_*=0$ (i.e. $\text{Im}i'_*\subset \text{kern}\, \ell'_*$). Indeed, consider an element
$f\in L^{0,n-1}_{\text{comp}}(\tilde{V},\,\mathcal{O}(D))\cap \text{kern}(\overline\partial)$ (where the sub-index $\text{comp}$ indicates that $f$ has compact support in $\tilde{V}$). Then $\ell'(f)_{|_{V_i}}$ is $\overline\partial$-closed, 
with compact support on $V_i$ and with ``polynomial blow-up''. By theorem 4.7, we know that there exists $v_i\in L^{0,n-1}_{2,\text{loc}}(V_i\setminus a_i)$, compactly supported in $V_i$ such that $\overline\partial v_i=\ell'(f)$ on $V_i\setminus a_i$. Hence $[\ell'f]=0$ in 
$^{\prime}H^{0,n-1}(X')$. 

\medskip
\noindent On the other hand, let $\ell'(g)=\overline\partial u$ with $u\in\,^{\prime}L^{0,n-2}_{(2)}(X')$, for some $g\in Z^{0,n-1}_{(2)}(\tilde{X},\,\mathcal{O}(D))$. Then $A:=g-\overline\partial \pi^*(\chi_0\,u)$ has compact support in $\tilde{V}$, so $[g]\in \text{Im}\,i'_*$. \;\;\;Q.E.D.

\medskip
\noindent{\bf{Remark 5.2.2:}} In the case of compact varieties $X$, we do not need to introduce the auxiliary spaces $^{\prime}L^{0,q}_{(2)}(X')$, ordinary local $L^2$-cohomology will do and Theorem 1.4 will be valid.  Moreover, in the case of projective surfaces we can prove the following corollary:

\begin{corollary} For projective surfaces  $X$ with finitely many isolated singularities, the map 

$$\tilde{T}: H^{0,1}_{(2)}(\tilde{X},\mathcal{O}(D))\to H^{0,1}_{(2)}(X')$$

\noindent of Theorem 1.4 is an isomorphism (the right-hand side $L^2$-cohomology is computed with respect to the restriction of the Fubini-Study metric in $X'$).  Here $\pi: \tilde{X}\to X$ is a desingularization of $X$ such that $E:=\pi^{-1}(\text{Sing}\,X)$ is a divisor 
with simple normal crossings, $Z:=\pi^{-1}(\text{Sing}\,X)$ is the unreduced exceptional divisor and $D:=Z-E$. 
\end{corollary}

\noindent{\it{Proof.}} We shall show that $H^1_E(\tilde{X},\mathcal{O}(D))=0$  by  showing that the map $H^1_E(\tilde{X},\,\mathcal{O}(D))\to H^1_E(\tilde{X},\,\mathcal{O}(Z))$ is injective.  Now $H^1_E(\tilde{X},\mathcal{O}(Z))\cong H^1_c(\tilde{U},\mathcal{O}(Z))$ where $\tilde{U}$ is a smooth strongly pseudoconvex neighborhood of $E$. The latter cohomology group is isomorphic to the dual of $H^1(\tilde{U},\mathcal{K}(-Z))$, which vanishes  by Takegoshi's or Silva's relative vanishing theorem, since $L_{-Z}$ is locally semi-positive with respect to $X$ (see example 11.22, page 56 in \cite{Dem1}, or \cite{Rup3} pages 24-25). Hence,  the proof of the corollary will be complete once we prove 

\begin{lemma} Under the assumptions of the corollary,  the map

$$H^1_E (\tilde{X},\,\mathcal{O}(D))\to H^1_E(\tilde{X},
\mathcal{O}(Z))$$ 

\noindent   is injective. 
\end{lemma}

\medskip
\noindent {\it{Proof.}} We introduce some auxiliary 1-cycles supported on $E=\cup_{i=1}^N  E_i$  and where $E_j$ are the irreducible components of $E$. For a special ordering of the irreducible
components of $E$ (to be determined later on), we set $D_0:=Z=\sum_{k=1}^N m_k\, E_k,\, D_j:=Z-\sum_{k=1}^j E_k$. Then $D_N=Z-E$. Consider the standard  short exact sequences of sheaves 

\begin{equation}\label{eq:finseslem}
0\to \mathcal{O}(D_j)\to \mathcal{O}(D_{j-1})\to \mathcal{O}_{E_j}(D_{j-1})\to 0.
\end{equation}

\medskip
\noindent Taking long exact sequence on cohomology with support on $E$ we obtain  for each $j\ge 1$

\begin{equation}\label{eq:finleslem}
...\to H^0_E(\tilde{X},\,\mathcal{O}_{E_j}(D_{j-1}))\to H^1_E(\tilde{X},\,\mathcal{O}(D_{j}))\to H^1_E(\tilde{X},\,\mathcal{O}(D_{j-1}))\to....
\end{equation}

\medskip
\noindent Suppose we were able to show that $E_j \,\cdot  D_{j-1}<0$\;\, for all $j\ge 1$ for some ordering of the irreducible components. Then $H^0_E(\tilde{X},\,\mathcal{O}_{E_j}(D_{j-1}))=H^0(E_j,\,\mathcal{O}_{E_j}(D_{j-1}))=0$.  This will imply  that each map $H^1_E(\tilde{X},\,\mathcal{O}(D_j))\to H^1_E(\tilde{X},\,\mathcal{O}(D_{j-1}))$ is injective for each $j=1,\,\cdots,  N$.  From this we can infer the injectivity of 
$H^1_E(\tilde{X},\,\mathcal{O}(D_N))\to H^1_E(\tilde{X},\,\mathcal{O}(D_0))$ which is precisely what we want in the lemma. 

\medskip
\noindent To conclude the proof of the lemma  it suffices to show that it is possible to rearrange the irreducible components $\{E_j\}$ of $E$ in such a way as to have $E_j \cdot D_{j-1}<0$ for all $j\ge 1$. The proof below is a generalization of the  proof of property a)  in the Appendix of \cite{PS1} (there they assumed that $E$ is connected, while we do not impose such a restriction).

\medskip\noindent Let $E_{(1)},\,\cdots, E_{(m)}$ denote the connected components of $E$.  We can write for each $1\le i\le m\; \;E_{(i)}:=\cup_{j\in J_i} E_j$  where $J_1,\,J_2,\cdots, J_m$ partition $\{1,2,\cdots, N\}$ and let $N_i:=|J_i|$.  

\medskip
\noindent  As the set $E$ is exceptional in $\tilde{X},$  let $\Phi: \tilde{X}\to Y$ be the blow-down map. By Proposition 4.6 in 
\cite{Lauf2}, since $\tilde{X}$ is normal, $Y$ is normal. But then, using Lemma 4.1 in \cite{Lauf2}, each connected component $E_{(i)}$ of $E$ ($1\le i\le m$) is mapped to a different point  $\{y_i\}$ of $Y$ under $\Phi$.  By theorem 4.4  in 
\cite{Lauf2}, the intersection matrix  for each connected component $E_{(i)}$ of $E$, denoted by  $S_{(i)}:=( E_{(i)}^j\cdot E_{(i)}^k)$\; for any ordering $E_{(i)}^1,\cdots, E_{(i)}^{n_i}$ of the irreducible components in  $E_{(i)},$  is negative definite. 

\medskip
\noindent Set $Z_{(i)}:=\sum_{k\in J_i} m_k\, E_{(i)}^k$. Let us observe that $E_{(i)}^j \cdot\, Z_ {(i)}=E^j_{(i)} \cdot\, Z$ for $j\in J_i$, since  irreducible components of $E$ that belong to different connected components do not intersect.  Following an idea of Gonzalez-Sprinberg (Lemma 2.1 in \cite{Gonz}), Pardon and Stern observed  (in the proof of property a) in the Appendix in \cite{PS1} as well as in Proposition 3.6 in \cite{PS3}) that for each irreducible component $E^k_{(i)}$ of $E_{(i)} $ one has $E^k_{(i)}\cdot Z_{(i)}\le 0$. Hence  we have $E^j_{(i)}\cdot Z_{(i)}=E^{j}_{(i)}\cdot Z\le 0$ for all 
$j\in J_i$. We claim now that there exists a  $j\in J_i$ such that $E^j_{(i)}\cdot Z_{(i)}<0$. Indeed, if for all $j\in J_i$ we had 
$E^j_{(i)}\cdot Z_{(i)} =0$ this would imply that for all $j\in J_i$ we have  $\sum_{k\in J_i} m_k E^{j}_{(i)} \cdot E^k_{(i)}=0$, which would contradict the 
negative definiteness of the matrix $S_{(i)}$. Hence there exists a $j\in J_i$ such that $E^{j}_{(i)} \cdot Z_{(i)}<0$. Let us call this  $E^j_{i}:=E_{i1}$.  Since $E_{(i)}$ is connected, we can inductively define $E_{i1},\,\cdots, E_{i\,N_i}$ such that $E_{ij}$ intersects some $E_{ik}$ for some $k<j;\,j>1$ and  such that

%there must exist a $k\in J_i$ with $k\neq j$ such that $E^j_{(i)}\cap E^k_{(i)}\neq \emptyset$ 
%(and hence $E^j_{(i)} \cdot E^k_{(i)}\ge 1$). Then $E^k_{(i)}\cdot (Z_{(i)}-E^j_{(i)})<0$. Let's call this $E^k_{(i)}:=E_{i2}$. This way %we can reorder/rename  the irreducible components  that belong in $E_{(i)}=\cup_{j=1}^{N_i} E_{ij}$ such that 

\begin{equation}\label{eq:choice}
E_{ij}\cdot (Z-\sum_{k=1}^{j-1} E_{ik})=E_{ij}\,\cdot (Z_{(i)}-\sum_{k=1}^{j-1} E_{ik}) <0.
\end{equation}

\medskip
\noindent Having ordered the irreducible components of E as above $\{E_{11},\,E_{12},\cdots E_{1\,N_1},\,E_{21},\,\cdots, E_{2\,N_2},\cdots E_{m1},\\
\noindent \cdots E_{m\,N_m}\}$ we can relabel them as $\{E_1,\,\cdots, E_{N_1}, E_{N_1+1}\,\cdots E_{N_1+N_2+1},\,\cdots E_{N_1+N_2+\cdots N_{m-1}+1},\,\cdots E_{N}\}$.  Then  we can show that  for all 
$j=1,\cdots N$ we have

$$ E_j \cdot (Z-\sum_{k<j} E_k)<0.$$\;\;\;
 
\noindent Indeed if $E_j:=E_{il}$ for some $i$ with $1\le i\le m$ and $1\le l\le N_i$, then taking into account that $E_{ik}\cdot Z=E_{ik}\cdot Z_{(i)}$ for all $k=1,\cdots N_i$,
 we can rewrite the above left-hand side as 

$$ E_{j}\cdot\,(Z-\sum_{k<j} E_k)=E_{il}\cdot Z_{(i)}-E_{il}\cdot (\underset{r<i,\,\\ 1\le \mu\le N_r} \sum E_{r\,\mu})-\sum_{k<l} E_{il}\cdot E_{ik}=E_{il}\cdot (Z_{(i)}-\sum_{k<l} E_{ik})<0,$$

\smallskip
\noindent where the last inequality follows from $(\ref{eq:choice})$ and the second term to the right hand side of the first equality vanishes due to the fact that 
irreducible components that belong to different connected components do not intersect and hence their intersection product is zero. \;\; Q.E.D.

\medskip\noindent 
{\bf{Remark 5.2.3}} Professor J\'anos Koll\'ar suggested an alternative proof of the vanishing of 
$H^1_E(\tilde{X},\,\mathcal{O}(D))$ in the case of a projective surface with a normal isolated singularity  based on duality and a 
strengthening of the Grauert-Riemenschneider vanishing theorem as it appears in Theorem 98, Chapter 2, page 51 in \cite{Kol}. For higher 
dimensional projective varieties with an  isolated singularity at a point $x$, Professor Koll\'ar reduced the vanishing of 
$H^{n-1}_E(\tilde{X},\,\mathcal{O}(D))$ to the vanishing of $H^{n-2}(E,\,\,{\mathcal{O}_{\tilde{X}}(Z)}_{\upharpoonright_{E}})$.  We present here 
a  proof that was inspired by his  argument:  

\smallskip
\noindent 
We consider  the short exact sequence 

$$
0\to \mathcal{O}_{\tilde{X}}(Z-E)\to \mathcal{O}_{\tilde{X}}(Z)\to \mathcal{O}_E(Z)\to 0
$$

\noindent  where $\mathcal{O}_E :=\mathcal{O}_{\tilde{X}}/\mathcal{O}_{\tilde{X}}(-E)$\;(sheaf supported on $E$)  and $\mathcal{O}_E(Z):=\mathcal{O}_E \otimes_{\mathcal{O}_{\tilde{X}}} \mathcal{O}_{\tilde{X}}(Z)$. 
Taking $\Gamma_E(\tilde{X},\,-)$ in the above short exact sequence we obtain a long exact sequence in cohomology with support in $E$

\begin{eqnarray}\label{eq:lese}
....\to H^{q-1}_E(\tilde{X},\,\mathcal{O}_{\tilde{X}}(Z)) \to H^{q-1}_E (\tilde{X},\,\mathcal{O}_E(Z))  \to H^{q}_E (\tilde{X},\,\mathcal{O}_{\tilde{X}}(Z-E))\to H^q_E(\tilde{X},\,\mathcal{O}_{\tilde{X}}(Z))\to...
\end{eqnarray}

\noindent By Karras' result we know that for $q<n$,  $H^{q}_E(\tilde{X},\,\mathcal{O}(Z))\cong H^q_c(\tilde{U},\,\mathcal{O}_{\tilde{X}}(Z))$, where $\tilde{U}$ is a smoothly bounded strongly pseudoconvex neighborhood of $E$ in $\tilde{X}$. The latter cohomology group (using Serre duality) is isomorphic to the dual of $H^{n-q}(\tilde{U},\,\mathcal{K}_{\tilde{X}}(-Z))$, which by Takegoshi's relative vanishing theorem will vanish if $n-q>0$. Hence for all $q<n$ we have $H^q_E(\tilde{X},\,\mathcal{O}(Z))=0$. 
We can then obtain from $(\ref{eq:lese})$ that 

$$
H^{n-1}_E(\tilde{X},\,\mathcal{O}(Z-E))\cong H^{n-2}(\tilde{X},\,\mathcal{O}_E(Z))\cong H^{n-2}(E,\,{\mathcal{O}_{\tilde{X}}(Z)}_{\upharpoonright_{E}}).
$$

\medskip\noindent  Hence 

$$
H^{n-1}_E(\tilde{X},\,\mathcal{O}(D))=0 \Longleftrightarrow H^{n-2}(E,\,{\mathcal{O}_{\tilde{X}}(Z)}_{\upharpoonright_{E}})=0.
$$

\smallskip
\noindent 

\medskip
\noindent {\bf{Remark 5.2.4}} For the local case now, exploiting the fact that a neighborhood of an isolated singularity embeds as an 
open subset in a projective variety, we can show that we always have $H^1_E(\tilde{X}_r,\,\mathcal{O}(D))=0$ in the $2$-dimensional case, 
and hence the map $T: H^{n-1}(\tilde{X}_r,\,\mathcal{O}(D))\to H^{0,n-1}_{(2)}(X'_r)$  of Theorem 1.2 is an isomorphism when $\text{dim}\,X=n=2$. 
It follows in an a similar way that the map $\tilde{T}$ of Theorem 1.4 is always an isomorphism when $n=2=\text{dim}\,X$.

\section{Proof of Corollary 1.6}

\medskip
\noindent In what follows, we use the assumptions and notation that were introduced  in the paragraph just above Corollary 1.6  (and in the paragraph above Theorems  1.3, 1.4)  in  section 1. 

\subsection{A description of the kernel of $i^n_*$} Recall that $i^n_*: H^n(\tilde{X},\,\mathcal{O}) \to H^{n}(\tilde{X},\,\mathcal{O}(D))$ is the map on cohomology induced by the sheaf inclusion $i:\mathcal{O}\to \mathcal{O}(D)$. Let us consider the following short exact sequence of sheaves

$$
0\to \mathcal{O}\to \mathcal{O}(D)\overset{\mu}\to \mathcal{O}_D(D)\to 0.
$$

\medskip
\noindent It yields two  long exact sequences on cohomology

\begin{equation*}
..H^{n-1}(\tilde{X},\,\mathcal{O})\overset{i^{n-1}_*}\to H^{n-1}(\tilde{X},\mathcal{O}(D))\overset{\mu^{n-1}_*}\to H^{n-1}(\tilde{X},\,\mathcal{O}_D(D))\overset{\delta}\to H^{n}(\tilde{X},\mathcal{O})\overset{i^n_*}\to H^n(\tilde{X},\mathcal{O}(D))\to H^{n}(\tilde{X},\mathcal{O}_D(D)) \\
\end{equation*}

\begin{equation*}
.. H^{n-1}(\tilde{U},\,\mathcal{O})\overset{^{\prime}  i^{n-1}_*}\to H^{n-1}(\tilde{U},\mathcal{O}(D))\overset{^{\prime} \mu^{n-1}_*}\to H^{n-1}(\tilde{U},\,\mathcal{O}_D(D))\overset{\delta}\to H^{n}(\tilde{U},\mathcal{O})=0\to..\phantom{ssasasasasdasdasdasdasa}\\
\end{equation*}

\noindent where $\tilde{U}=\pi^{-1}(U)$ with $U$ a disjoint  union of  smoothly bounded  strongly pseudoconvex neighborhoods of the singular points $\{a_j\}_{j=1}^m$. The vanishing of $H^n (\tilde{U},\mathcal{O})$ is due to Siu's theorem in \cite{Siu}.  Also, $H^n(\tilde{X},\,\mathcal{O}_D(D))=H^n(|D|,\,\mathcal{O}_D(D))=0$, since the support of D, denoted by  $|D|,$  is an $(n-1)$-dimensional variety.  From the exactness of the first long exact sequence we know that $\text{kern}(i^n_*)=\text{Im}(\delta)$,\;and that $\text{Im}(\mu^{n-1}_*)  \overset{\Psi}  \cong \dfrac{H^{n-1}(\tilde{X},\mathcal{O}(D))}{i^{n-1}_*(\,H^{n-1}(\tilde{X},\mathcal{O})\,)}.$  We also have the following short exact sequence that defines $\text{kern}\,i^n_*$

$$
0\to \text{Im}(\mu^{n-1}_*) \overset{I} \to H^{n-1}(|D|,\,\mathcal{O}_D(D)) \overset{\delta}\to \text{kern}\,i^n_*\to 0,
$$

\noindent where $I$ is the inclusion map. 

\medskip
\noindent From the exactness of the second long exact sequence above we obtain that 

$$
H^{n-1}(\tilde{U},\mathcal{O}_D(D))=H^{n-1}(|D|,\mathcal{O}_D(D)) \overset{\Omega}\cong \frac{H^{n-1}(\tilde{U},\mathcal{O}(D))}{\text{kern}(^{\prime} \mu^{n-1}_*)}= \frac{H^{n-1}(\tilde{U},\,\mathcal{O}(D))}{^{\prime} i^{n-1}_* (H^{n-1}(\tilde{U},\,\mathcal{O}))}.
$$

\medskip
\noindent  The short exact sequence that defines $\text{kern}(i^n_*)$ can be rewritten now as follows:

$$
0\to 
{\dfrac{H^{n-1}(\tilde{X},\mathcal{O}(D))}{i^{n-1}_*(H^{n-1}(\tilde{X},\mathcal{O}))}}\overset{\nu}\to { \dfrac{H^{n-1}(\tilde{U},\mathcal{O}(D))}{^{\prime}i^{n-1}_*(H^{n-1}(\tilde{U},\mathcal{O}))}}\overset{\delta\circ \Omega^{-1}}\longrightarrow \text{kern}(i^n_*)\to 0,
$$

\noindent for some injective map $\nu:=\Omega\circ I\circ \Psi^{-1}$. 

\smallskip
\noindent 
Using the commutativity of the right grid in the very first diagram of section 4.5 and Karras' results, we can conclude that the maps $i^{n-1}_*,\,^{\prime} i^{n-1}_*$ (denoted in section 4.5 as $m_*$)  are  injective. By identifying $i^{n-1}_*(H^{n-1}(\tilde{X},\mathcal{O}))$ with $H^{n-1}(\tilde{X},\mathcal{O})$\, and $ ^{\prime} i^{n-1}_*( H^{n-1}(\tilde{U},\mathcal{O}))$ with $H^{n-1}(\tilde{U},\mathcal{O})$ from the above short exact sequence we see that if 

\medskip
\noindent  $$\text{dim}_{\mathbb{C}} { \dfrac{H^{n-1}(\tilde{U},\mathcal{O}(D))}{H^{n-1}(\tilde{U},\mathcal{O})}}\neq \text{dim}_{\mathbb{C}} 
{\dfrac{H^{n-1}(\tilde{X},\mathcal{O}(D))}{H^{n-1}(\tilde{X},\mathcal{O})}},$$

\noindent then the kernel of $i^n_*$ would be non-trivial. 

\subsection{Proof of Corollary 1.6} 

\medskip
\noindent We need to show that  i)\; $\text{kern}(\phi^n_*) \subset \text{kern}(i^n_*)$ and   ii)\; $\text{kern}(i^n_*)\subset \text{kern}(\phi^n_*)$.  In what follows we shall think of $i^n_* :  H^{n}(\tilde{X},\mathcal{O})\to H^n(\tilde{X},\mathcal{O}(D))\cong H^{0,n}_{(2)}(\tilde{X}, L_D)$ 
as the map on cohomology induced by the sheaf map $\mathcal{O}\to \mathcal{O}(L_D)$, sending $f\to f\otimes s$, where $s$  is the canonical section of $L_D$ introduced in section 2.2.

\medskip
\noindent To prove i),  let $[c]\in H^{n}(\tilde{X},\mathcal{O})$ such that $\phi^n_* ([c])=[0]$. Without loss of generality we can assume that $[c]$ can be represented by an element $g\in L^{0,n}_{(2)}(\tilde{X})$ with $g=0$ in $\tilde{U}$ (since we can solve $\overline\partial t=g$ in a neighborhood of  $\overline{\tilde{U}}$  we can replace $g$ by $g-\overline\partial(\xi\,t^o)$, where $\xi$ is a 
cut-off function with $\xi=1$ on $\overline{\tilde{U}}$ and $t^o$ denotes trivial extension by zero outside $\tilde{U}$).  Now, 
$\phi^n(g)=\overline\partial u$ for some $u\in L^{0,n-1}_{(2)}(X')$. Using a cut-off function $\chi \in C^{\infty}(X)$  with $\chi=1$ near 
the singular locus $A$ and $\text{supp}\,\chi\subset U$, we can rewrite $\phi^n(g)$ in $U'$ as

$$(*)\;\;\;\;\;\;\;\;\;\;\;\;\;\;\;\;\;\;\;\;\;\phi^n(g)=\overline\partial u=\overline\partial \chi \wedge u+
\overline\partial \left((1-\chi)\,u\right).
$$

\medskip
\noindent This is possible since  $g$ was taken to be $0$ on $\tilde{U}$  yielding   $u\in Z^{0,n-1}_{(2)}(U')$. Using 
the surjectivity of the map $\ell_*$ on $\text{Im}\,j_*$ (Lemma 4.9 in our paper), we know that there exists an $A\in Z^{0,n-1}_{(2)}(\tilde{U},\,L_D)$ and a $v\in L^{0,\,n-2}_{2,\text{loc}}(\tilde{U})$ such that $u_{\upharpoonright U}=(\pi^{-1})^* (A \cdot s^{-1})+\overline\partial v$. Setting $\tilde{\chi}=\chi\circ \pi$  and applying  $\pi^*$ on both sides  of $(*),$ \; we obtain on $\tilde{U}\setminus E$   

\begin{eqnarray*}
(**) \;\;\;\;\;\;\;\;\;\;\;\;\;\;\;g&=&\overline\partial \tilde{\chi}\wedge \pi^* u+\overline\partial \left((1-\tilde{\chi})\,\pi^* u \right)\\ 
&=& \overline\partial \tilde{\chi}\wedge \left( A\cdot s^{-1}+\overline\partial \pi^* v \right)+\overline\partial \left
((1-\tilde{\chi})\,\pi^* u \right).
\end{eqnarray*}

\medskip
\noindent  From $(**)$\; we obtain that $g\otimes s=\overline\partial \left( \tilde{\chi}\,A^o-\overline\partial \tilde{\chi}\wedge \pi^* v^o \otimes s+(1-\tilde{\chi})\,\pi^*u \otimes s\right)=\overline\partial B$, with $B=\tilde{\chi}\,A^o -\overline\partial \tilde{\chi}\wedge \pi^* v^o\otimes s+(1-\tilde{\chi})\,\pi^*u\otimes s\in L^{0,n-1}_{(2)}(\tilde{X},\L_D)$, since $\pi$ is a quasi-isometry from $\{\tilde{\chi}<1\}$ onto $\{\chi<1\}$.  Hence $i^n_*([g])=[0]$. 

\medskip
\noindent To prove ii) we consider again an element $[c]\in \text{kern}(i^n_*)$. As in the proof of i), we can assume that this class may be represented by an element $g\in L^{0,n}_{(2)}(\tilde{X})$ with $g=0$ in  $\tilde{U}$. By assumption we have that there 
exists an element $A\in L^{0,n-1}_{(2)}(\tilde{X}, \,L_D)$ such that $g\otimes s=\overline\partial A=\overline\partial \tilde{\chi}\,\wedge A+\overline\partial \left(\,(1-\tilde{\chi})\,A\right)$ and thus we have on $\tilde{U}\setminus E$

$$
(***)\phantom{dfsdfsdsdfsd}g=\overline\partial \tilde{\chi}\wedge A\cdot s^{-1}+\overline\partial \left((1-\tilde{\chi})\, (A\cdot s^{-1})\right)
$$

\medskip
\noindent since  $g=0$ on $\tilde{U}$  and therefore  $A\in Z^{0,n-1}_{(2)}(\tilde{U},\,L_D)$. Now by the surjectivity of $\ell_*$ on $\text{Im}\,j_*$ (Lemma 4.9 in our paper) we know that there exist elements $t\in Z^{0,n-1}_{(2)}(U')$ and $v\in L^{0,n-2}_{2,\,\text{loc}}(U')$ such that on $U'$ we have  $(\pi^{-1})^*( A\cdot  s^{-1})=t+\overline\partial v$.  Applying $\phi^n$ on\; 
$(***)$\; we can express  $\phi^n(g)$ on $X'$ as $\phi^n(g)=\overline\partial \chi \wedge (t^o+\overline\partial v^o)+\overline\partial \left((1-\chi) \phi^n(A\cdot s^{-1})\right)=\overline\partial \left(\chi\,t^o-\overline\partial \chi \wedge v^o-(1-\chi)\, \phi^n(A\cdot s^{-1})\right)=\overline\partial C$ where $C:=\chi\,t^o-\overline\partial \chi \wedge v^0-(1-\chi)\, \phi^n(A\cdot s^{-1})\in L^{0,n-1}_{(2)}(X')$. Hence $\phi^n_*([g])=[0]$. \;\;\;\;Q.E.D.

\medskip
\noindent Using Corollary1.6, we can easily show the existence of the following commutative diagram

$$
\xymatrix{ 
&H^{n} (\tilde{X},\,\mathcal{O})\ar[d]^{i^n_*} \ar[r]^{\phi^*_n}& H^{0,n}_{(2)}(X')\\
&H^{n}(\tilde{X},\,\mathcal{O}(D)) \ar@{-->}[ur]\\
}
$$

\noindent  As the maps $\phi^n_*, \,i^n_*$ are surjective and $\text{kern}(i^n_*)=\text{kern}(\phi^n_*)$  the dotted map will be  an isomorphism.  Thus $H^{0,n}_{(2)}(X')\cong H^n(\tilde{X},\,\mathcal{O}(D))$. \;\;\;Q.E.D.

\medskip
\section{Examples}
\noindent The purpose of this section is to produce various
examples for which we have or not vanishing of
$L^2$-$\overline\partial$-cohomology groups.

\medskip
\noindent In \cite{FOV} we showed that whenever $0$ was a
Cohen-Macaulay point of a pure $n$-dimensional complex analytic
variety with $n\ge 3$, then we have $H^{0,q}_{(2)}(X'_r)=0$ for
all $q$ with $1\le q\le n-2$ (this result was obtained using
Theorem  2.4  and an extension theorem of cohomology classes by
Scheja). Classical examples of Cohen-Macaulay singularities are
rational singularities of dimension $n\ge 2$ (Corollary 4.3 in
\cite{Kar}-attributed to Kempf). Recall that in a complex space
$X$, a normal point $p\in X$  is called {\it{rational}} if given a
resolution of singularities $\pi:\tilde{X}\to X$  we have that
$\left(R^{i} \pi_*\mathcal{O}_{\tilde{X}}\right)_p=0$ for all
$i>0$. It follows from Hironaka's work that the condition on
$R^i\pi_*\mathcal{O}_{\tilde{X}}$ is independent of the choice of
$\tilde{X}$. There is a plethora  of rational singularities as the
following examples suggest.

\medskip
\noindent\;{\it{Example 1: Quotient singularities.}} We have the
following theorem proven by Burns:

\begin{theorem}(Proposition 4.1 in \cite{Bur}) Let $M$ be a complex manifold and $G$ a properly
discontinuous group of automorphisms of $M$. Then $X=M/G$ has only
rational singularities.
\end{theorem}

\noindent This implies all the double point singularities in
dimension 2 ($A_k, D_k, E_6, E_7, E_8$) are rational singularities.

\medskip
\noindent\;{\it{Example 2: Arnold's singularities (example 3.4 in
\cite{Bur}).}} These are direct generalizations of the rational
double points of dimension 2.

\medskip
\noindent {\it{Example 3: Some affine cones over smooth projective
hypersurfaces (Example 1.2  in \cite{Bur}).}}  Let $V\subset \mathbb{CP}^n$ be a smooth
hypersurface of degree $d$, described by the nonsingular homogeneous
polynomial $F(Z)$ in the homogeneous coordinates $Z_1,\cdots Z_{n+1}$.
Let $X$ be defined by $F(Z)=0$ in $\mathbb{C}^{n+1}$. $X$ is called
the affine cone over $V$. Let $\mathcal{O}_V(-1)$ denote the
restriction of the universal line bundle on $V,\;\,p: \mathcal{O}_V(-1)\to V$ the corresponding projection map and let $\tilde{X}$
denote the total space of $\mathcal{O}_V(-1)$. Now, the map $\pi: \tilde{X} \to X$ is the
contraction of the zero section of $\mathcal{O}_V(-1)$, hence $\pi$
a resolution of singularities of $X$ with $\pi^{-1}(0)=V$. In the algebraic category, using
this, we can see that $0\in X$ is a rational singularity if and only
if $d\le n$. Indeed, by the Leray spectral sequence of $\pi$ we have
that $H^j (\tilde{X},\,\mathcal{O}_{\tilde{X}})=H^0(X,\,R^j\pi_*
\mathcal{O}_{\tilde{X}})=(R^i\,\pi_*\,\mathcal{O}_{\tilde{X}})_0$.  Using the Leray spectral sequence for $p$ (the 
fact that $R^i\,p_*\mathcal{O}_{\tilde{X}}=0$ for all $i>0,$   hence $H^j (\tilde{X},\,\mathcal{O}_{\tilde{X}})\cong H^j (V,\,p_* \mathcal{O}_{\tilde{X}})$ )\footnotemark\footnotetext{This vanishing is a consequence of  Exercises 8.1-8.2 page 252 in \cite{Har}.}  
and expanding cohomology classes into Taylor series along the fibers  
of $\mathcal{O}_V(-1)$ we get that 
$H^j(\tilde{X},\,\mathcal{O}_{\tilde{X}})=\oplus_{k\ge 0}
\,H^j(V,\,\mathcal{O}_V(k))$. Using the fact that the canonical line
bundle of $V$ is given by  $K_V=[(d-n-1)H_{|_V}]$ (adjunction
formula) where $H$ is the hyperplane bundle in $\mathbb{CP}^n$,
Serre's duality and the dual version of Kodaira's vanishing theorem
for negative line bundles, we see that all these cohomology groups
vanish for $j>0$ if $d\le n$. In this case $0$ is a rational
singularity for the affine cone over $V$.

\medskip
\noindent On the other hand, when $d=n+1$ the  calculation in
Example 3 will yield
$H^{n-1}(\tilde{X},\,\mathcal{O})=H^{n-1}(V,\,\mathcal{O})$. If
$\text{dim}_{\mathbb{C}} H^{n-1}(V,\,\mathcal{O})\neq 0$ then one
produces a non-rational singularity. This happens for example if $V$
is any Riemann surface in $\mathbb{CP}^2$, of genus $g\ge 1$.

\medskip
\noindent {\bf{Remark:}} The local $L^2$-$\overline\partial$-
cohomology groups are completely determined in the case of affine
cones over smooth projective varieties. In this case $Z=|Z|$ and Theorem 1.2 in our paper or 
Theorem 7.1 of Ruppenthal in \cite{Rup3} guarantees that
$H^{0,q}_{(2)}(X'_r)\cong H^q(\tilde{X}_r,\,\mathcal{O})$ for all
$q$ with $1\le q\le n-1$.

\medskip \noindent \subsection{Some remarks on non Cohen-Macaulay
spaces} There are many irreducible complex analytic spaces that
are not Cohen-Macaulay. In a very interesting paper \cite{SV}
St\"uckrad and Vogel constructed a wealth of examples of smooth
projective varieties $V$ (Proposition 9 in \cite{SV}) whose affine
cone over $V$, denoted by $X(V)$ and abbreviated by $X$ when there
is no confusion, had the property that its local ring at the
vertex of the cone was not Cohen-Macaulay. The precise
construction is as follows:

\medskip
\noindent Let $d\ge 3$ and consider the variety $W\subset
\mathbb{CP}^{d-1}$ defined by the equation $z_0^d+z_1^d+\cdots
z_{d-1}^d=0$. Let $V$ be the Segr\'e embedding of
$W\times\mathbb{CP}^1$ in $\mathbb{CP}^{2d-1}$. In \cite{SV} it is
shown that the local ring $\mathcal{O}_{X,\,0}$ ($0$ is the vertex
of the affine cone $X$ over $V$) is a normal non Cohen-Macaulay
ring. Andreatta and Silva used this construction to produce in
\cite{AS} another example of non-rational singularity.

\medskip
\noindent Using a K\"unneth formula for Segr\'e products one has

$$
H^{d-1}(V,\,\mathcal{O}_V(k))\cong \oplus_{r+s=d-1}
\left(H^r(W,\,\mathcal{O}_W(k))\otimes
H^s(\mathbb{CP}^1,\,\mathcal{O}_{\mathbb{CP}^1}(k))\right)
$$

\noindent By the version of Kodaira's vanishing theorem for negative
line bundles and using also the fact that $K_W\cong \mathcal{O}_W$
by the adjunction formula, the dual of the  cohomology groups
$H^r(W,\,\mathcal{O}_W(k))$ can be computed

$$(H^r(W,\,\mathcal{O}_W(k)))'\cong H^{d-2-r}(W, K_W\otimes
\mathcal{O}_W(-k))\cong H^{d-2-r}(W,\mathcal{O}_W(-k))=0$$

\noindent for $r>0$ and $k>0$. By the K\"unneth formula above we
have

$$
H^{d-1}(V,\,\mathcal{O}_V(k))\cong
H^{d-2}(W,\,\mathcal{O}_{W}(k))\otimes
H^1(\mathbb{CP}^1,\,\mathcal{O}_{\mathbb{CP}^1}(k))=0.
$$

\noindent since
$H^1(\mathbb{CP}^1,\,\mathcal{O}_{\mathbb{CP}^1}(k))=0$ for all
$k\ge -1$. We know that there must exist an $i$ with $1\le i\le
d-2,$ such that $H^i(\tilde{X},\mathcal{O}_{\tilde{X}})\neq 0$,
where $\tilde{X}$ is a desingularization of the affine cone $X$
over $V$ (otherwise, the vertex of the cone would have been a
rational singularity, thus Cohen-Macaulay). As mentioned in the
remark above, using Theorem 1.2   we have that
$H^{0,d-1}_{(2)}(X'_r)=0$ while, for some $i$ with $1\le i\le d-2$
we have that $H^{0,i}_{(2)}(X'_r)\neq 0$ (using Theorem 1.1 ). 

\medskip
\noindent The vertex of the cone in this example is a new type of
singularity called {\it{weakly rational}} singularity. Recall that
in an $n$-dimensional complex space $(X,\,\mathcal{O}_X)$, a point
$p\in X$ is called a {\it{weakly rational singularity}} of $X$ if
$\left(R^{n-1} \pi_*\mathcal{O}_X\right)_p=0$. As before
$\pi:\tilde{X}\to X$ is a resolution of $X$. The above definition
is independent of the resolution $\pi$. Andreatta and Silva and
Yau studied to what extent Laufer's results on rational
singularities from \cite{Lauf1} generalize to this category of
singularities. It is clear that when $n=2$, then the definitions
of weakly rational and rational coincide. For higher dimensional
singularities this is no longer true as the following example
shows:

\medskip
\noindent {\it{Example.}} Consider a compact Riemann surface $V$
of genus  $1$ and $F$ a sufficiently
negative\footnotemark\footnotetext{Let us recall Proposition 2.1
from \cite{Lauf3}: Let $N$ be a vector bundle over a compact
Riemann surface $A$. Suppose that $N=\oplus_{i=1}^{n-1} L_i,$
where $L_i$ is a line bundle of Chern class $c_i$ on $A$. Then $A$
is exceptional in $N\,\Longleftrightarrow\,c_i<0,\;\,1\le i\le n-1$.}
vector bundle of rank  $r\ge 2$  over $V$. Let $\tilde{X}$ denote
the total space of $F$. Let $\pi: \tilde{X}\to X$ be the
blow-down of  $V$ in $\tilde{X}$ and $x:=\pi(V)$. Then, $x$ is
weakly rational since $H^{\text{dim}_x
X-1}(\tilde{X},\,\mathcal{O}_{\tilde{X}})=0$ but not rational
since $H^1(\tilde{X},\,\mathcal{O}_{\tilde{X}})\neq 0$. To show
this we need the following facts: a) Proposition 26 from
Andreotti-Grauert \cite{AG}-that discusses filtrations of
cohomology groups of vector bundles over complex manifolds and
their associated graded complex- and asserts that $\text{Grad}\,
H^i(\tilde{X},\,\mathcal{O}_{\tilde{X}})\cong
\oplus_{k=0}^{\infty} H^i(V,\,\mathcal{O}(F^*)^k),$
\;b)\;\; the fact that $\text{dim}_{\mathbb{C}}
H^1(V,\,\mathcal{O}(F^*)^k)=1$ for $k=0$ and $0$
otherwise, and last \;c)\;the knowledge that the cohomology groups
$H^q(V,\,\mathcal{O}(F^*)^k)=0$ for all $k\ge 0$ and
$q\ge 2$.

\medskip
\noindent Moreover $\mathcal{O}_{X,\,x}$ is not Cohen-Macaulay;
recall that if the homological codimension
$\text{codh}\,\mathcal{O}_{X,\,x}=\text{dim}_x X\ge 3\ge 2+1$ then
by Theorem 3.1, page 37 in \cite{BS} we should have
$H^i_{x}(X,\,\mathcal{O})=0$ for all $i\le 2$. From the local
cohomology exact sequence and taking into account that $X$ is
Stein, we see that $H^2_x(X,\,\mathcal{O})=H^1(X\setminus
\{x\},\,\mathcal{O})\cong H^1(\tilde{X}\setminus
V,\,\mathcal{O})$. Using Karras' results the latter cohomology
group is isomorphic to $H^1(\tilde{X},\,\mathcal{O})$. The latter
space is nonzero from earlier computations. Hence
$\text{codh}\,\mathcal{O}_{X,\,x}\neq \text{dim}_x \,X$.

\medskip \noindent Using Theorem 1.1 and these
calculations we obtain $\text{dim} H^{0,1}_{(2)}(U\setminus
\{x\})\neq 0$, where $U$ is a small Stein neighborhood of $x$ with
smooth boundary and $H^{0,q}_{(2)}(U\setminus \{x\})=0$ for $2\le
q\le \text{dim}_x \,X-2$ (if $\text{dim}_x X\ge 4$). Now by further blowing up $V$ inside $\tilde{X}$ (for example blowing-up the ideal sheaf of $V$) we can obtain a  manifold $\hat{X}$, a map $p:\hat{X}\to \tilde{X}$ such that $p^{-1}(V)=E$ ($E$ is non-singular
as it is locally isomorphic to $V\times \mathbb{P}^{s-1}$, with $s=\text{codim}(V,\,\tilde{X})$ and locally principal--see for example Theorem 8.24, page 186 in \cite{Har}). Then we obtain a resolution map $\hat{p}: \hat{X}\to X$, such that $\hat{p}^{-1}(x)=E$. Using Theorem 1.2   (as $-D$ is locally semi-positive with respect to $X$), we
obtain $H^{0,\text{dim}_x X-1}_{(2)}(U\setminus\{x\})\cong 
H^{\text{dim}_x X-1}(\hat{p}^{-1}(U),\,\mathcal{O}(D))$ where
$D=Z-|Z|$ and $Z=m\,E$ is the unreduced divisor $\hat{p}^{-1}(x)$ and $|Z|=E$.  If $m=1$ then $\mathcal{O}(D)\cong \mathcal{O}$, hence $H^{\text{dim}_x X-1} (\hat{X},\,\mathcal{O}((m-1)\,E))=H^{\text{dim}_x X-1}(\hat{X},\,\mathcal{O})=0$ by the calculation above. 
Determining the vanishing or not of $H^{\text{dim}_x X-1} (\hat{X},\,\mathcal{O}((m-1) E))$  when  $m\ge 2$ is slightly more involved. 

\medskip
\noindent
We know that $\mathcal{O}_{\hat{X}}(E)_{\upharpoonright_E}=N_{E|\hat{X}}=\mathcal{O}_{P(N_{V| \tilde{X}})}(-1)$ (see Proposition 12.4 in section 12, Chapter VII of \cite{Dem2}), where $N_{V|\tilde{X}}$ is the normal bundle of $V$ in $\tilde{X}$ and $\mathcal{O}_{P(N_{V|\tilde{X}})}(-1)$ is the tautological line bundle over $E=P(N_{V|\tilde{X}})$.  Here the projectivized normal bundle $P(N_{V|\tilde{X}})$ is defined by considering lines in $N_{V|\tilde{X}}$.

\medskip
\noindent
Using observation  $\beta)$  in  the Characterization of Exceptional Sets (section 3.1)  on page 10,   and   Proposition 26 from \cite{AG}, we have

\begin{eqnarray*}
 H^{\text{dim}_x X-1}(\hat{p}^{-1}(U),\,\mathcal{O}(D))&\cong&   H^{\text{dim}_x X-1} (\hat{X},\,\mathcal{O}((m-1)\,E)),\\
\text{Grad}\, H^{\text{dim}_x X-1} (\hat{X},\,\mathcal{O}((m-1)\,E)) &=&\oplus_{k\ge 0}\;\;H^{\text{dim}_x X-1} (E,\,\mathcal{O}((m-1)\,E)_{\upharpoonright_E}\,\otimes \mathcal{O}_E(1)^k\,),\\
\end{eqnarray*}

\medskip\noindent where $\mathcal{O}_E(1)=\mathcal{O}_E(-1)^*$ (the dual of the tautological line bundle over $E$). Hence, 

\begin{eqnarray}
\text{Grad}\, H^{\text{dim}_x X-1}(\hat{X},\,\mathcal{O}((m-1)\,E))&=& \oplus_{k\ge 0} H^{\text{dim}_x X-1}(E,\,\mathcal{O}_E\,(k-m+1))\\
&=& \oplus_{k\ge 0}\, H^0\,(E,\,K_E\otimes \mathcal{O}_E\,(m-1-k))
\end{eqnarray}
 
\medskip
\noindent  The first equality follows from the fact that $\mathcal{O}((m-1)\,E)_{\upharpoonright E}=\mathcal{O}_E(-(m-1))$. The second follows from Serre duality. 
For a fixed rank $r$ of the vector bundle $F$, the vanishing or not of  $H^{\text{dim}_x X-1} (\hat{X},\,\mathcal{O}((m-1)\,E))$  depends on the multiplicity $m$ of the divisor $Z$.  To determine what happens for  $m\ge 2$  we need to recall some facts about projectivized vector bundles.   We summarize them in the following Proposition:

\begin{proposition} Let $\mathcal{W}$ be a holomorphic vector bundle of rank $r\ge 2$ over a complex manifold $M$, let $P(\mathcal{W}):=\frac{\mathcal{W}\setminus 0}{\mathbb{C}^*}$ be the projectivized bundle  and let $p: P(\mathcal{W})\to M$ be the associated projection map. Let 
$\mathcal{L}^{-1}_{P(\mathcal{W})}:=\{((z,\,[v]),\,\lambda\,v)\;|\;\;(z,\,[v])\in P(W),\;\lambda\in \mathbb{C}\}$ be the tautological line bundle associated to $\mathcal{W}$  (it is a sub-bundle of $p^*\,\mathcal{W}$) , let $\mathcal{L}_{P(\mathcal{W})}$ denote its dual bundle (``hyperplane bundle'') and let $\odot^k \mathcal{W}^*$ denote the $k$-th symmetric product of $\mathcal{W}^*$.  Then we have: 

\medskip
\noindent i) (Proposition 2.2 in \cite{KO}\footnotemark\footnotetext{Kobayashi-Ochiai's $L(\mathcal{W})$ corresponds in our notation to $\mathcal{L}_{P(\mathcal{W})}^{-1}$.} or Theorem 3.5 in \cite{CW})  The canonical bundle $K_{P(\mathcal{W})}$ is given by 

$$ K_{P(\mathcal{W})}=p^*\,(\text{det}\,\mathcal{W}^*\otimes K_M)\otimes \mathcal{L}_{P(\mathcal{W})}^{-r},$$

\noindent where $\mathcal{L}^{-r}_{P(\mathcal{W})}$ is the dual of the $r$-fold tensor product of the ``hyperplane bundle'' $\mathcal{L}_{P(\mathcal{W})}$. 

\medskip
\noindent ii) (Lemma 3.1 in \cite{Har2} \footnotemark\footnotetext{Hartshorne's $\bf{P}(\mathcal{W})$ corresponds in our notation to 
$P(\mathcal{W}^*)$.}  or Theorem attributed to Grothendieck \footnotemark\footnotetext{La theorie des classes de Chern, Bull. S. M. F, 
tome 86, (1958), 137-154.}  on page 403 in \cite{CW})  For all $m\ge 0,\;i>0$  and for any coherent analytic sheaf $\mathcal{S}$ on $M$,  we have

\medskip

\medskip
\noindent
ii.a) \;\;\;\;\;\;\;\;\;\;\;\;\;\;\;\;\;\;\;\;\;\;\;\;\;\;\;\;\;\;\;\;\;\;\;\;\;\;\;\;\;\;\;\;\;$  p_*(\mathcal{L}_{P(\mathcal{W})}^m\otimes p^*\,\mathcal{S})\cong \odot^m \mathcal{W}^*\otimes S,$

\medskip
\noindent
ii.b) \;\;\; \;\;\;\;\;\;\;\;\;\;\;\;\;\;\;\;\;\;\;\;\;\;\;\;\;\;\;\;\;\;\;\;\;\;\;\;\;\;\;\;\;\;\;$ R^i\,p_*(\mathcal{L}_{P(\mathcal{W})}^m \otimes p^*\,\mathcal{S})=0.$ 

\medskip
%\noindent iii) (Theorem, page 30 in the preprint version of \cite{CW})  For all $m>0,\;i\neq r-1$, we have 

%$$ R^i\,(\mathcal{L}_{P(\mathcal{W})}^{-m}\otimes p^*\,\mathcal{S})=0$$

%\noindent while for $i=r-1$, we have 

%$$ R^{r-1}\,p_*(\mathcal{L}_{P(\mathcal{W})}^{-m}\otimes p^*\mathcal{S})\cong \odot^{m-r}\,\mathcal{W}\otimes \mathcal{S},$$

%\noindent where by convention $\odot^k \mathcal{W}=0$, if $k>0$. 

\end{proposition}

\medskip
\noindent  In what follows we will abbreviate the normal bundle $N_{V|\tilde{X}}$ of $V$ in $\tilde{X}$ by $NV$.  Recall that $E=P(NV)$ and $p: P(NV)\to V$  is the associated projection map.  Identifying $\mathcal{O}_E(-1)=\mathcal{L}_{P(NV)}^{-1}$ and using Proposition 7.2 i) (taking into account that $K_V=0$),  we can rewrite the right-hand side of (23)  as

$$\oplus_{k\ge 0} H^0(E,\,p^*(\text{det}( NV)^*)\otimes \mathcal{L}_{P(NV)}^{-r-k+m-1}).$$

\medskip
\noindent Now, if $-r-k+m-1= 0$ (this can happen for example for certain values of $k$ when $m\ge 3$ and $r=2$ or more generally when $m\ge r+1$), we obtain from Proposition  7.2 ii.a)  

\begin{equation}\label{eq:van}
H^0(E,\,p^*(\text{det} (NV)^*)\otimes \mathcal{L}_{P(NV)}^{-r-k-1+m})\cong H^0(V,\, \text{det}\,(NV)^*\otimes \mathcal{O}_V),
\end{equation}

\noindent since $\odot^0 (NV)^*=\mathcal{O}_V$. Now  $\text{det}(NV)^*$ is a line bundle of positive degree over $V$ so  by Riemann-Roch we know it has a non-zero section. Hence, the  right-hand side of $(\ref{eq:van})$ (and therefore $H^{\text{dim}_x X-1}(\hat{X},\,\mathcal{O}((m-1)\,E)))$   does not vanish.  

\medskip
\noindent On the other hand,  if $-r-k-1+m<0$ \; (this can happen for all non-negative integers $k$,  when $m=2$ and $r\ge 2$ or more generally when $m<r+1$) then taking into account that $p_* (\mathcal{L}_{P(E)}^{-\nu}\otimes p^* \mathcal{S})=0$ for all $\nu\ge 1$ and all coherent analytic 
sheaves $\mathcal{S}$ on $V$,  we obtain that  $H^0(E,\,p^*(\text{det}\,(NV)^*)\otimes \mathcal{L}_{P(NV)}^{-r-k-1+m})=0$. Hence $H^{\text{dim}_x X-1}(\hat{X},\,\mathcal{O}((m-1)\,E))=0$   in this case (i.e.  when $m<r+1$).

%\begin{comment}

\end{document}